\documentclass[twocolumn]{autart}    

\usepackage{graphicx}          
\usepackage[tbtags]{amsmath}
\usepackage[misc]{ifsym}
\usepackage{mathrsfs}
\usepackage{amssymb}
\usepackage{amsmath}
\usepackage{amsfonts}
\usepackage{graphicx}
\usepackage{enumerate}
\usepackage{setspace}
\usepackage{color,xcolor}

\def\IR{\text{\rm l\negthinspace R}}
\def\IE{\text{\rm l\negthinspace E}}
\def\IP{\text{\rm l\negthinspace P}}

\begin{document}

\begin{frontmatter}

\title{An optimal control problem for mean-field forward-backward stochastic differential equation with noisy observation\thanksref{footnoteinfo}} 

\thanks[footnoteinfo]{This work is supported in part by 
the NSF of China under Grants 11371228, 61422305, 61304130 and 61573217, by the NSF for Distinguished Young Scholars of Shandong
Province of China under Grant JQ201418, and by the Research Fund for the
Taishan Scholar Project of Shandong Province of China. The material in this work was partially presented at the 35th Chinese Control Conference, July 27-29, 2016, Chengdu. Corresponding author: Hua Xiao.}

\author[Add1]{Guangchen Wang}\ead{wguangchen@sdu.edu.cn},               
\author[Add1,Add2]{Hua Xiao}\ead{xiao\_hua@sdu.edu.cn},
\author[Add1]{Guojing Xing}\ead{xgjsdu@sdu.edu.cn}  

\address[Add1]{School of Control Science and Engineering, Shandong University, Jinan 250061, China}             
\address[Add2]{School of Mathematics and Statistics, Shandong University, Weihai 264209, China}                  


\begin{abstract}                          
This article is concerned with an optimal control problem derived by mean-field forward-backward stochastic differential equation with noisy observation, where the drift coefficients of the state equation and the observation equation are linear with respect to the state and its expectation. The control problem is different from the existing literature on optimal control for mean-field stochastic systems, and has more applications in mathematical finance, e.g., asset-liability management problem with recursive utility, systematic risk model. Using a backward separation method with a decomposition technique, two optimality conditions along with two coupled forward-backward optimal filters are derived. Several linear-quadratic optimal control problems for mean-field forward-backward stochastic differential equations are studied. Closed-form optimal solutions are explicitly obtained in detailed situations.
\end{abstract}

\end{frontmatter}

\section{Introduction}

\subsection{Notation}

We denote by $T>0$ a fixed time horizon, by $\IR^m$ the $m$-dimensional Euclidean space, by $|\cdot|$ (resp.
$\langle\cdot,\cdot\rangle$) the norm (resp. scalar product) in a Euclidean space, by $A^\top$ (resp. $A^{-1}$) the transposition
(resp. reverse) of $A$, by $S^m$ the set of symmetric $m\times m$ matrices with real elements, by $f_x$ the partial derivative of
$f$ with respect to $x$, and by $C$ a positive constant, which can be different from line to line. Let $(\Omega, \mathscr F, (\mathscr F_t)_{0\leq t\leq T}, \IP)$ be a complete filtered probability space, on which are given an $\mathscr F_t$-adapted standard Brownian motion $(w_t, \tilde w_t)$ with values in $\IR^{r+\tilde r}$ and a Gaussian random variable $\xi$ with mean $\mu_0$ and covariance matrix $\sigma_0$. $(w, \tilde w)$ is independent of $\xi$. If $A\in S^m$ is positive (semi) definite, we write $A>(\geq)0$. If $x:[0, T]\rightarrow\IR^m$ is uniformly bounded, we write $x\in\mathscr L^\infty(0,T; \IR^m)$. If $x: \Omega\rightarrow
\IR^m$ is an $\mathscr F_T$-measurable, square-integrable random variable, we write $x\in \mathscr L^2_{\mathscr F}(\IR^m)$. If $x: [0, T]\times\Omega\rightarrow \IR^m$ is an $\mathscr F_t$-adapted, square-integrable process, we write $x\in \mathscr L^2_{\mathscr
F}(0, T; \IR^m)$. We also adopt similar notations for other processes, Euclidean spaces and filtrations.

\subsection{Motivation}

Now consider an asset-liability management problem of a firm. Let the dimension $n=k=r=\tilde r=1.$ Denote by $\IE$ the expectation with respect to $\IP$, by $v_t$ the control strategy of the firm, by $x^v_t$ the cash-balance, and by $\bar l^v_t$ the liability process. Norberg \cite{nr} described the liability process by a Brownian motion with drift. The model, however, is not just the one we want. In fact, it is possible that the control strategy and the mean of the cash-balance can influence the liability process, due to the complexity of the financial market and the risk aversion behavior of the firm. Such an example can be found in Huang et al. \cite{hww1}, where the liability process depends on a control strategy (for example, capital injection or withdrawal) of the firm. Along this line, we proceed to improve the liability process here. Suppose that $\bar l^v_t$ satisfies a linear stochastic differential equation (SDE, in short) of the form
$$-d\bar l^v_t=(\bar a_t\IE x^v_t+b_tv_t+\bar b_t)dt+c_tdw_t.$$
Here $\bar a$, $b$, $\bar b$, $c$, $a$, $f$, $g$ and $h$ are deterministic and uniformly bounded. $\bar b_t$ and $c_t$ denote the liability rate and the volatility coefficient, respectively. Suppose that the firm owns an initial investment $\xi$, and only invests in a money account with the compounded interest rate $a_t$. Then the cash-balance of the firm is
$$x^v_t=e^{\int_0^ta_sds}\left(\xi-\int_0^te^{-\int_0^sa_rdr}d\bar l^v_s\right).$$
It follows from It\^{o}'s formula that
$$
\left\{\begin{aligned}
dx^v_t=&\ \left(a_tx^v_t+\bar a_t\IE x^v_t+b_tv_t+\bar b_t\right)dt+c_tdw_t,\\
x^v_0=&\ \xi.
\end{aligned}
\right.
$$
Note that, if $b_t=1$, $\bar b_t=0$, $a_t=-\bar a_t=\hbox{const.}$ and $c_t=\hbox{const.}$, then the cash-balance equation is just the systematic risk model of inter-bank borrowing and lending introduced in Carmona et al. \cite{cfs}. Besides the systematic risk model, the equation can also be reduced to an air conditioning control model in energy-efficient buildings. See, e.g., Example 2 in Djehiche et al. \cite{dtt} for more details.

Due to the discreteness of account information, it is possible for the firm to partially observe the cash-balance by the stock price
$$
\left\{\begin{aligned}
dS^v_t=&\ S^v_t\left[\left(f_tx^v_t+g_t+\frac{1}{2}h^2_t\right)dt+h_td\tilde w_t\right],\\
S^v_0=&\ 1.
\end{aligned}
\right.
$$
Set $Y^v_t=\log S^v_t$. It holds that $Y^v$ is governed by
$$
\left\{\begin{aligned}
dY^v_t=&\ \left(f_tx^v_t+g_t\right)dt+h_td\tilde w_t,\\
Y^{v}_0=&\ 0.
\end{aligned}
\right.
$$
Suppose that the firm has triple performance objectives. The first two ones are to minimize the total cost of $v$ over $[0, T]$ and to minimize the risk of $x^v_T$. Assume that the risk is measured by $\IE\left[(x^v_T-\IE x^v_T)^2\right].$ The third one is to maximize the utility $y^v_t$ resulting from $v$. Without loss of generality, define
$$
y^v_t=\IE\left[x^v_T+\int_t^TG(s, y^v_s, v_s)ds\bigg|\mathscr F_t\right],
$$
where $G$ is Lipschitiz continuous with respect to $(y, v)$, and $G(s, 0, 0)\in\mathscr L^2_{\mathscr F}(0, T; \IR)$ for $0\leq s\leq T$. We emphasize that the current utility $y^v_t$ depends not only on the instantaneous control $v_t$, but also on the future utility and control $(y^v_s, v_s)$, $t\leq s\leq T$. This shows the difference between the utility $y^v$ and the standard additive utility, and hence, $y^v$ is called as a stochastic differential recursive utility in Duffie and Epstein \cite{de}. Then the asset-liability management problem with recursive utility is stated as follows.

\textbf{Problem (AL).} Find a $\sigma\{Y^v_s; 0\leq s\leq t\}$-adapted and square-integrable process $v_t$ such that
$$J[v]=\frac{1}{2}\IE\left[\int_0^TB_tv^2_tdt+H(x^v_T-\IE x^v_T)^2-2Ny^v_0\right]$$ is minimized. Here $B>0$ and $B^{-1}$ are deterministic and uniformly bounded. $H$ and $N$ are non-negative constants. $y^v_0$ is the value of $y^v_t$ at time $0$.

Let us now turn to the recursive utility $y^v_t$ again. According to El Karoui et al. \cite{epq}, $y^v_t$ admits the backward stochastic differential equation (BSDE, in short)
$$
\left\{\begin{aligned}
-dy^v_t=&\ G(t, y^v_t, v_t)dt-z^v_tdw_t-\tilde z^v_td\tilde w_t,\\
y^v_T=&\ x^v_T.
\end{aligned}
\right.
$$
With the BSDE, Problem (AL) can be rewritten as an optimal control problem derived by forward-backward stochastic differential equation (FBSDE, in short) with noisy observation. It is possible to work out one more asset-liability management problem. We omit the details to limit the length of this article.

\subsection{Problem statement}

Motivated by the examples, we study an optimal control problem for FBSDE with noisy observation. Consider a controlled FBSDE
$$
\left\{\begin{aligned}
dx^v_t=&\ \left(a_tx^v_t+\bar a_t\IE x^v_t+b(t, v_t)\right)dt+c_tdw_t,\\
-dy^v_t=&\ (\alpha_tx^v_t+\bar\alpha_t\IE x^v_t+\beta_ty^v_t+\bar\beta_t\IE y^v_t+\gamma_tz^v_t\\
&\ +\bar\gamma_t\IE z^v_t+\tilde\gamma_t\tilde z^v_t+\bar{\tilde\gamma}_t\IE\tilde z^v_t+\psi(t, v_t))dt\\
&\ -z^v_tdw_t-\tilde z^v_td\tilde w_t,\\
x^v_0=&\ \xi, \quad y^v_T=\rho x^v_T+\bar\rho\ \IE x^v_T,
\end{aligned}
\right.
$$
where $(x^v, y^v, z^v, \tilde{z}^v)$ is the sate, $v$ is the control, and $(w, \tilde w)$ is the Brownian motion. Since the mean of the state influences
the state equation, we call the equation a mean-field FBSDE, or a McKean-Vlasov FBSDE. Assume that $(x^v, y^v, z^v, \tilde{z}^v)$ is partially observed through
$$
\left\{\begin{aligned}
dY^v_t=&\ \left(f_tx^v_t+\bar f_t\IE x^v_t+g(t, v_t)\right)dt+h_td\tilde w_t,\\
Y^{v}_0=&\ 0.
\end{aligned}
\right.
$$
The cost functional is
$$J[v]=\IE\left[\int_0^Tl(t, x^v_t, \IE x^v_t, v_t)dt+\phi(x^v_T, \IE x^v_T)+\varphi(y^{v}_{0})\right].$$
Here $v_t$ is required to be $\sigma\{Y^v_s; 0\leq s\leq t\}$-adapted and to satisfy $\IE\sup_{0\leq t\leq T}|v_t|^2<+\infty.$ $a$, $\bar a$, $b$, $c$, $\alpha$, $\bar \alpha$, $\beta$, $\bar \beta$, $\gamma$, $\bar \gamma$, $\tilde{\gamma}$, $\bar{\tilde{\gamma}}$, $\psi$, $\rho$, $\bar \rho$, $f$, $\bar f$, $g$, $h$, $l$, $\phi$ and $\varphi$ will be specified in Section 2. Our problem is to select an admissible control $v$ to minimize $J[v]$. We denote the problem by Problem (MFC), where ``MF" and ``C\," are the capital initials of ``mean-field" and ``control", respectively.

To solve Problem (MFC), it is natural to use dynamic programming and maximum principle. The dynamic programming, however, does not hold even if the BSDE and the observation equation are not present, mainly due to the inclusion of the mean of the state, which leads to the time inconsistency. We instead study the maximum principle for optimality of Problem (MFC).

\subsection{Briefly historical retrospect and contribution of this paper}

Mean-field theory provides an effective tool for investigating the collective behaviors arising from individuals' mutual interactions in various different fields, say, finance, game, engineering. Since the independent introduction by Lasry and Lions \cite{ljlp} and Huang et al. \cite{hcm1,hcm2}, the mean-field theory has attracted more and more attention. Let us now briefly recall some latest developments which are related to Problem (MFC).

Although the study of mean-field SDE has a long history with the pioneering works of Kac \cite{km} and McKean \cite{mhp}, mean-field type control is a rather new research direction. In 2001, Ahmed and Ding \cite{ad} used the Nisio nonlinear operator semigroup to obtain an extended dynamic programming. By dual techniques, maximum principles for several mean-field SDEs with full information were derived. See, e.g., Buckdahn et al. \cite{bdl}, Li \cite{lj}, Hafayed and Abbas \cite{ha}, Shen et al. \cite{sms}, Djehiche et al. \cite{dtt}. Subsequently, Meyer-Brandis et al. \cite{moz}, Hafayed et al. \cite{haa1,haa2} studied the partial information case, where noisy observation and filtering are excluded. As applications of the derived maximum principles, \cite{moz,wwz1,hnz} solved mean-variance problems with full and partial information. Yong \cite{yj1} studied a linear-quadratic (LQ, in short) optimal control problem for mean-field SDE with full information. Further, Yong \cite{yj2} investigated the time-inconsistency feature of the LQ problem, and obtained both open-loop and closed-loop equilibrium solutions. Later, Huang et al. \cite{hly} extended the LQ problem to the case of infinite horizon. For the discrete-time counterpart of the LQ problem, please refer to Elliott et al. \cite{eln}, Ni et al. \cite{nel,nzl} and the references therein for more details. It is worth pointing out that the investigation of mean-field type control is also partially inspired by the interest in mean-field game. If we only focus on a single decision maker, also called a representative agent, mean-field game can be regarded as mean-field type control. Generally speaking, an exact Nash equilibrium for mean-filed game with a large number of decision makers is rarely available except for special cases (see, e.g., Carmona et al. \cite{cfs}). It is highly desirable to find a good approximation of this Nash equilibrium. Please refer to Carmona et al. \cite{cdl}, Tembine et al. \cite{tzb}, Bensoussan et al. \cite{bfy}, etc. for more details on different types of approximation equilibrium. See also Bensoussan et al. \cite{bsyy} for a comprehensive study of a general LQ mean field game.

Both mean-field type control and mean-field game lead to mean-field FBSDE. Buckdahn et al. \cite{bdlp} studied the well-posedness of a decoupled mean-field FBSDE using a limit approach. Bensoussan et al. \cite{byz}, Carmona and Delarue \cite{cd} extended \cite{bdlp} to the fully coupled mean-field FBSDE case in terms of a continuation method introduced in Peng and Wu \cite{pw}. Mean-field FBSDE is a well-defined dynamic system, it is very natural and appealing to study control and game problems for mean-field FBSDEs as well as their applications. To our knowledge, there is only a few literature on this topic. For example, Li and Liu \cite{lrlb} studied an optimal control problem for fully coupled mean-field FBSDE. Hafayed et al. \cite{htb} obtained a maximum principle for mean-field FBSDE with jump. Huang et al. \cite{hww2} studied an LQ game with a linear mean-field BSDE system and a quadratic cost functional. \cite{htb,hww2} also provided some applications in mean-variance and recursive utility problems.

In this paper, we are interested in studying Problem (MFC). Compared with the above literature, this problem has several new features as follows.

\begin{itemize}

\item The state $(x^v, y^v, z^v, \tilde z^v)$ satisfies a mean-field FBSDE rather than a mean-field SDE, and is only partially observed by a noisy process. This endows Problem (MFC) more practical meanings in reality.

\item Unlike those control models solved in Bensoussan \cite{ba2}, the classical separation principle does not work here, mainly due to the fact that the mean square error of filtering of BSDE depends on the control in general.

\item The state equation involves the mean of the state, and thus, Problem (MFC) can not be studied by transforming it into a standard control problem for FBSDE. This feature can be supported by Example 2.2 in this paper.
\end{itemize}

There is a few papers related to Problem (MFC). Let us make a brief comment on them. Wang et al. \cite{wzz} posted a partially observable mean-field type optimal control problem for SDE. They used a backward separation method and a probability transformation to decouple a circular dependence between the control and the observation first, and then derived a necessary condition for optimality. The result was further generalized in Wang et al. \cite{wwz1} by the backward separation method with an approximation technique. Later, Hu et al. \cite{hnz} studied an optimal control problem for mean-field SDE with jump. Zhang \cite{zh} addressed the case with correlated state and observation noises.  We emphasize that the approach applied in \cite{hnz,wwz1,wzz,zh} is based on at least one of the assumptions below.

\begin{itemize}
\item The state satisfies an SDE rather than an FBSDE.

\item The drift term of the observation equation is uniformly bounded with respect to the state and the control.

\item The control has no effect on the observation.

\item The control $v$ satisfies $\IE\sup_{0\leq t\leq T}|v_t|^\ell<+\infty$, $\forall\ell>0$.
\end{itemize}
Clearly, Problem (MFC) does not meet these assumptions. Another approach is desired to develop to address Problem (MFC). In \cite{wwx2}, Wang et al. studied an LQ control problem for classical FBSDE (i.e., the dynamics of the FBSDE does not depend on the probability distribution of the state). Inspired by Bensoussan \cite{ba1}, they solved the LQ problem by combining a decomposition technique with the backward separation method. Recently, our further study on the approach provided in Wang et al. \cite{wwx2} finds out the availability of the approach to some nonlinear control problems with noisy observations, say, Problem (MFC). In this paper, we will show how to use the approach to address Problem (MFC). See also Wang et al. \cite{wwx1} for other developments about partially observable optimal control for FBSDE.

The rest of this article is organized as follows. In Section 2, we carefully formulate Problem (MFC) first, and then provide illustrative examples and preliminary results. In Section 3, we obtain two optimality conditions and two coupled forward-backward optimal filtering equations. In Section 4, we study an LQ case of Problem (MFC) and obtain a feedback representation of optimal control. In Section 5, we explicitly solve an asset-liability management problem with noisy observation, and work out an illustrative numerical example. Some concluding remarks and proofs of the preliminary results are given in Section 6 and Appendix, respectively.

\section{Problem formulation and preliminary}

One difficulty to study Problem (MFC) is there is a circular dependence between the control $v$ and the observation $Y^v$, which results in the unavailability of classical variation.
Here we will adopt a decomposition technique, similar to those of \cite{ba1,wwx2}, to overcome the difficulty.
Define $(x^0, y^{0}, z^{0}, \tilde{z}^{0})$ and $Y^0$ by
\begin{equation}
\label{se1} \left\{\begin{aligned}
dx^0_t=&\ \left(a_tx^0_t+\bar a_t\IE x^0_t\right)dt+c_tdw_t,\\
-dy^0_t=&\ \left(\alpha_tx^0_t+\bar\alpha_t\IE x^0_t+\beta_ty^0_t+\bar\beta_t\IE y^0_t+\gamma_tz^0_t\right.\\
&\ \left.+\bar\gamma_t\IE z^0_t+\tilde\gamma_t\tilde z^0_t+\bar{\tilde\gamma}_t\IE\tilde z^0_t\right)dt\\
&\ -z^0_tdw_t-\tilde z^0_td\tilde w_t,\\
x^0_0=&\ \xi, \quad y^v_T=\rho x^0_T+\bar\rho\ \IE x^0_T,
\end{aligned}
\right.
\end{equation}
and
\begin{equation}
\label{oe1} \left\{\begin{aligned}
dY^0_t=&\ \left(f_tx^0_t+\bar f_t\IE x^0_t\right)dt+h_td\tilde w_t,\\
Y^{0}_0=&\ 0,
\end{aligned}
\right.
\end{equation}
where
$a$, $\bar a\in\mathscr L^\infty(0, T; \IR^{n\times n})$,
$c\in\mathscr L^\infty(0, T; \IR^{n\times r})$,  $\alpha$, $\bar \alpha\in\mathscr L^\infty(0, T; \IR^{m\times n})$, $\beta$, $\bar \beta$ $\in\mathscr L^\infty(0, T;$ $\IR^{m\times m})$, $f,$ $\bar f\in\mathscr L^\infty(0, T; \IR^{\tilde r\times n})$, $h, h^{-1}\in\mathscr L^\infty(0, T; \IR^{\tilde r\times\tilde r})$,
$\gamma=(\gamma_{1},\cdots,\gamma_{r})$, $\bar{\gamma}=(\bar{\gamma}_{1},\cdots,\bar{\gamma}_{r})$,
$\tilde{\gamma}=(\tilde{\gamma}_{1},\cdots,\tilde{\gamma}_{\tilde{r}})$, $\bar{\tilde{\gamma}}=(\bar{\tilde{\gamma}}_{1},\cdots,\bar{\tilde{\gamma}}_{\tilde{r}})$,
$z^{0}=(z^{0}_{1},\cdots,z^{0}_{r})$,
$\tilde{z}^{0}=(\tilde{z}^{0}_{1},\cdots, \tilde{z}^{0}_{\tilde{r}})$,
$\gamma_{j}$, $\bar{\gamma}_{j}$, $\tilde{\gamma}_{j}$,  $\bar{\tilde{\gamma}}_{j}$ $\in\mathscr L^\infty(0, T;$ $\IR^{m\times m})$ for $j=1, \cdots, r$ or $\tilde{r}$, and $\rho$, $\bar{\rho}\in\IR^{m\times n}$ are constant matrices. Here we use the simplified notation $\gamma_tz^0_t \triangleq
\sum_{j=1}^{r}\gamma_{j t}z^{0}_{j t}$. Similarly, it is also applicable for the notations $\tilde{\gamma}_t\tilde{z}^0_t $, $\bar{\gamma}_t \IE {z}^0_t $, $\bar{\tilde{\gamma}}_t \IE {\tilde{z}}^0_t, \cdots $.

Let $v\in\mathscr L^2_{\mathscr F}(0, T; \IR^k)$ be a control process. Define $(x^{v, 1}, y^{v, 1}, z^{v, 1}, \tilde{z}^{v, 1})$ and $Y^{v, 1}$  by
\begin{equation}
\label{se2} \left\{\begin{aligned}
\dot x^{v, 1}_t=&\ a_tx^{v, 1}_t+\bar a_t\IE x^{v, 1}_t+b(t, v_t),\\
-dy^{v, 1}_t=&\ \left(\alpha_tx^{v, 1}_t+\bar\alpha_t\IE x^{v, 1}_t+\beta_ty^{v, 1}_t+\bar\beta_t\IE y^{v, 1}_t\right.\\
&\ +\gamma_tz^{v, 1}_t+\bar\gamma_t\IE z^{v, 1}_t+\tilde\gamma_t\tilde z^{v, 1}_t+\bar{\tilde\gamma}_t\IE\tilde z^{v, 1}_t\\
&\ \left. +\psi(t, v_t)\right)dt-z^{v, 1}_tdw_t-\tilde z^{v, 1}_td\tilde w_t,\\
x^{v, 1}_0=&\ \xi, \quad y^{v, 1}_T=\rho\, x^{v, 1}_T+\bar\rho\ \IE x^{v, 1}_T,
\end{aligned}
\right.
\end{equation}
and
\begin{equation}
\label{oe2} \left\{\begin{aligned}
\dot Y^{v, 1}_t=&\ f_tx^{v, 1}_t+\bar f_t\IE x^{v, 1}_t+g(t, v_t),\\
Y^{v, 1}_0=&\ 0,
\end{aligned}
\right.
\end{equation}
where $g: [0,T]\times\IR^k\rightarrow\IR^{\tilde r}$ satisfies $\IE\int_0^T|g(t, v_t)|^2dt<+\infty$, $b: [0,
T]\times\IR^k\rightarrow\IR^n$ and $\psi: [0,
T]\times\IR^k\rightarrow\IR^m$ are continuous and continuously differentiable with
respect to $t$ and $v$, respectively, $b_v\in\mathscr L^\infty(0, T; \IR^{n\times k})$ and  $\psi_v\in\mathscr L^\infty(0, T; \IR^{m\times k})$. Since (\ref{se1}) and (\ref{se2}) are decoupled, it is easy to see from Buckdahn et al. \cite{bdlp} that (\ref{se1}), (\ref{oe1}), (\ref{se2}) and (\ref{oe2}) have unique solutions, respectively. Define
\begin{equation}
\label{so}
\begin{aligned}
x^v_t=&\ x^0_t+x^{v,1}_t, \quad  y^v_t=y^0_t+y^{v,1}_t,\\
z^v_t=&\ z^0_t+z^{v,1}_t, \quad  \tilde{z}^v_t=\tilde{z}^0_t+\tilde{z}^{v,1}_t
\end{aligned}
\end{equation}
and
\begin{equation}\label{eq1}
 Y^v_t=Y^0_t+Y^{v,1}_t.
\end{equation}
It follows from It\^o's formula that $(x^v, y^v, z^v, \tilde{z}^v)$ and $Y^v$ uniquely solve
\begin{equation}\label{se}
\left\{\begin{aligned}
dx^v_t=&\ \left(a_tx^v_t+\bar a_t\IE x^v_t+b(t, v_t)\right)dt+c_tdw_t,\\
-dy^v_t=&\ \left(\alpha_tx^v_t+\bar\alpha_t\IE x^v_t+\beta_ty^v_t+\bar\beta_t\IE y^v_t+\gamma_tz^v_t\right.\\
&\ \left. +\bar\gamma_t\IE z^v_t+\tilde\gamma_t\tilde z^v_t+\bar{\tilde\gamma}_t\IE\tilde z^v_t+\psi(t, v_t)\right)dt\\
&\ -z^v_tdw_t-\tilde z^v_td\tilde w_t,\\
x^v_0=&\ \xi, \quad y^v_T=\rho\, x^v_T+\bar\rho\, \IE x^v_T,
\end{aligned}
\right.
\end{equation}
and
\begin{equation}
\label{oe} \left\{\begin{aligned}
dY^v_t=&\ \left(f_tx^v_t+\bar f_t\IE x^v_t+g(t, v_t)\right)dt+h_td\tilde w_t,\\
Y^{v}_0=&\ 0,
\end{aligned}
\right.
\end{equation}
respectively. Let $\mathscr F^{Y^v}_t=\sigma\{Y^v_s; 0\leq s\leq t\} \ \hbox{and}\ \mathscr F^{Y^0}_t=\sigma\{Y^0_s; 0\leq s\leq t\}.$
Note that the first observable filtration depends on the control $v$. However, the second one is not the case. We now give a definition of admissible control. Let $U$ be a nonempty convex subset of $\IR^k$, and let $\mathscr U^0_{ad}$ be the collection of all
$\mathscr F^{Y^0}_t$-adapted processes with values in $U$
such that $\IE\sup_{0\leq t\leq T}|v_t|^2<+\infty.$

\textbf{Definition 2.1.} A control $v$ is called admissible, if
$v\in\mathscr U^0_{ad}$ is $\mathscr F^{Y^v}_t$-adapted. The set of
all admissible controls is denoted by $\mathscr U_{ad}$.

With the definition, it follows from the equality (\ref{eq1}) that

\textbf{Proposition 2.1.} \emph{For any $v\in\mathscr U_{ad}$,
$\mathscr F^{Y^v}_t=\mathscr F^{Y^0}_t.$}

It implies that the control $v$ has no effect on the observation $Y^v$, i.e., the circular dependence between $v$ and $Y^v$
is decoupled.

The cost functional is in the form of
\begin{equation}
\label{cf}
\begin{aligned}
J[v]=&\ \IE\left[\int_0^Tl(t, x^v_t, \IE x^v_t, v_t)dt\right.\\
&\ \left. +\phi(x^v_T, \IE x^v_T)+\varphi(y^{v}_{0})\right],
\end{aligned}
\end{equation}
where $l: [0, T]\times\IR^{n+n}\times U\rightarrow\IR$,
$\phi: \IR^{n+n}\rightarrow\IR$ and $\varphi: \IR^{m}\rightarrow\IR$ are continuously differentiable with
respect to $(x, \bar x, v)$, $(x, \bar x)$ and $y$, respectively, and there is a constant $C>0$ such that
$$
\begin{aligned}
|\phi(x, \bar x)|\leq&\ C(1+|x|^2+|\bar x|^2),\\
|\phi_x(x, \bar x)|+|\phi_{\bar x}(x, \bar x)|\leq&\  C(1+|x|+|\bar x|),\\
|l(t, x, \bar x, v)|\leq&\  C(1+|x|^2+|\bar x|^2+|v|^2),\\
|l_\chi(t, x, \bar x, v)|\leq&\  C(1+|x|+|\bar x|+|v|),\\
|\varphi(y)|\leq&\ C(1+|y|^2),\\
|\varphi_{y}(y)|\leq&\ C(1+|y|^2)
\end{aligned}
$$
with $\chi=x, \bar x, v$. Then the optimal control problem for mean-field FBSDE is restated as follows.

\textbf{Problem (MFC).} Find a $u\in\mathscr U_{ad}$ such that
$J[u]=\inf_{v\in\mathscr U_{ad}}J[v]$ subject to (\ref{se}),
(\ref{oe}) and (\ref{cf}). Any $u$ satisfying the equality is called an optimal control of Problem (MFC), and $(x^u, y^u, z^u, \tilde{z}^u)$ and $J[u]$ are called
the optimal state and the optimal cost functional, respectively.

Note that the above decomposition technique is restricted to special structures of state and observation equations, say, the case that (\ref{se}) and (\ref{oe}) are linear with respect to $(x^v, y^v, z^v, \tilde z^v)$, the diffusion coefficient of (\ref{se}) is deterministic, and the drift coefficient of (\ref{oe}) is independent of $(y^v, z^v, \tilde z^v)$. It is worth investigating the availability of the technique to decompose more general state and observation equations in the future.

Next, let us show more new features of Problem (MFC) by two simple examples. Roughly speaking, Example 2.1 tells us that
Problem (MFC) is possibly applied to solve a partially observable optimal control problem for mean-field SDE with stochastic coefficients in certain situations; Example 2.2 reveals that Problem (MFC) is not a trivial extension to a partially observable optimal control problem for FBSDE without mean-field term.

\textbf{Example 2.1.} Let $\alpha_t=\bar\alpha_t=\beta_t=\bar\beta_t=\bar\gamma_t=\bar{\tilde\gamma}_t=0$ in Problem (MFC). Then (\ref{se}) is reduced to
\begin{equation}
\label{rse}\left\{\begin{aligned}
dx^v_t=&\ \left(a_tx^v_t+\bar a_t\IE x^v_t+b(t, v_t)\right)dt+c_tdw_t,\\
-dy^v_t=&\ (\gamma_tz^v_t+\tilde\gamma_t\tilde z^v_t+\psi(t, v_t))dt-z^v_tdw_t-\tilde z^v_td\tilde w_t,\\
x^v_0=&\ \xi, \quad y^v_T=\rho x^v_T+\bar\rho\IE x^v_T.
\end{aligned}
\right.
\end{equation}
Solving the BSDE in (\ref{rse}), we get
\begin{equation}
\label{ys} y^v_0=\IE\left[\langle\eta_T, \rho x^v_T+\bar\rho\IE x^v_T\rangle
+\int_0^T\langle\eta_t, \psi(t, v_t)\rangle dt\right]
\end{equation}
with
$$
\left\{\begin{aligned}
d\eta_t=&\ \gamma_t\eta_tdw_t+\tilde\gamma_t\eta_td\tilde w_t,\\
\eta_0=&\ I_m,
\end{aligned}
\right.
$$
where $I_m$ is an $m$-dimensional vector with all components being $1$. Plugging (\ref{ys}) into (\ref{cf}), we have
$$
\begin{aligned}
\mathcal J[v]=&\ \IE\left[\int_0^Tl(t, x^v_t, \IE x^v_t, v_t)dt+\phi(x^v_T, \IE x^v_T)\right.\\
&\ +\varphi\left(\IE\left[\langle\eta_T, \rho x^v_T+\bar\rho\IE x^v_T\rangle\right.\right.\\
&\ \left.\left.\left. +\int_0^T\langle\eta_t, \psi(t, v_t)\rangle dt\right]\right)\right].
\end{aligned}
$$
Then Problem (MFC) is reduced to minimize $\mathcal J[v]$ over $\mathscr U_{ad}$ subject to (\ref{oe}) and the SDE in (\ref{rse}).
It is worth noting that we start with a control model with deterministic coefficients, but we end up with a control model with stochastic coefficients. The interesting phenomena is caused by the introduction of the BSDE in (\ref{rse}). Just because of this, maybe it provides a potential method to investigate a control problem for mean-field SDE with stochastic coefficients under certain conditions, i.e., we can change it into an equivalent control problem for mean-field FBSDE with deterministic coefficients. The details of how to make use of this potential method will be shown in our future publications, because they beyond the scope of the present paper.

\textbf{Example 2.2.} Find a $v\in\mathscr U_{ad}$ such that
$$
\begin{aligned}
J[v]=&\ \IE\left[\int_0^T\left(\langle A_tx^v_t, x^v_t\rangle+\langle \bar A_t\IE x^v_t, \IE x^v_t\rangle\right.\right.\\
&\ \left.+\langle B_tv_t, v_t\rangle\right)dt+\langle Hx^v_T, x^v_T\rangle\\
&\ \left.+\langle \bar H\IE x^v_T, \IE x^v_T\rangle+\langle My^v_0, y^v_0\rangle\right]
\end{aligned}
$$
is minimized, subject to (\ref{se}) and (\ref{oe}) with the assumption
$b(t, v_t)=\tilde{b}_tv_t+\bar b_t$, $\psi(t, v_t)=\tilde{\psi}_tv_t+\bar\psi_t$, $g(t, v_t)=\tilde{g}_tv_t+\bar g_t$,
where $A, \bar A\in\mathscr
L^\infty(0, T; S^n)$, $B\in\mathscr L^\infty(0, T; S^k)$, $H, \bar
H\in S^n$, $M\in S^m$, $A, B, H, M\geq0$, $A+\bar A, H+\bar H \geq0$, $\tilde{b}\in\mathscr L^\infty(0, T; \IR^{n\times k})$,
$\bar b\in\mathscr L^\infty(0, T; \IR^n)$, $\tilde{\psi}\in\mathscr L^\infty(0, T; \IR^{m\times k})$,
$\bar \psi\in\mathscr L^\infty(0, T; \IR^m)$, $\tilde{g}\in\mathscr L^\infty(0, T; \IR^{\tilde r\times k})$ and $\bar g\in\mathscr L^\infty(0, T; \IR^{\tilde r})$.
For simplicity, we denote the LQ problem by \textbf{Example (MFLQ)}.

Taking expectations on both sides of (\ref{se}) and (\ref{oe}), we have
$$
\left\{\begin{aligned}
\dot{\IE x^v_t}=&\ \frac{d}{dt}\IE x^v_t=(a_t+\bar a_t)\IE x^v_t+\tilde{b}_t\IE v_t+\bar b_t,\\
-\dot{\IE y^v_t}=&\ -\frac{d}{dt}\IE y^v_t=(\alpha_t+\bar\alpha_t)\IE x^v_t+(\beta_t+\bar\beta_t)\IE y^v_t\\
                 &\ +(\gamma_t+\bar\gamma_t)\IE z^v_t+(\tilde\gamma_t+\bar{\tilde\gamma}_t)\IE\tilde z^v_t\\
                 &\ +\tilde{\psi}_{t}\IE v_t+\bar\psi_{t},\\
\IE x^v_0=&\ \mu_0,  \quad \IE y^v_T=(\rho +\bar\rho)\ \IE x^v_T,
\end{aligned}
\right.
$$
and
$$\left\{\begin{aligned}
\dot{\IE Y^v_t}=&\ \frac{d}{dt}\IE Y^v_t= (f_t+\bar f_t)\IE x^v_t+\tilde{g}_t\IE v_t+\bar g_t,\\
\IE Y^{v}_0=&\ 0,
\end{aligned}
\right.
$$
respectively. Then
$$
\left\{\begin{aligned}
d(x^v_t-\IE x^v_t)=&\ \left[a_t(x^v_t-\IE x^v_t)+\tilde{b}_t(v_t-\IE v_t)\right]dt\\
&\ +c_tdw_t,\\
-d(y^v_t-\IE y^v_t)=&\ \left[\alpha_t(x^v_t-\IE x^v_t)+\beta_t(y^v_t-\IE y^v_t)\right.\\
&\ \left.+\gamma_t(z^v_t-\IE z^v_t)+\tilde\gamma_t(\tilde z^v_t-\IE\tilde z^v_t)\right.\\
&\ \left.+\tilde{\psi}_t(v_t-\IE v_t)\right]dt\\
&\ -z^v_tdw_t-\tilde z^v_td\tilde w_t,\\
x^v_0-\IE x^v_0=&\ \xi-\mu_0, \\
y^v_T-\IE y^v_T=&\ \rho(x^v_T-\IE x^v_T),
\end{aligned}
\right.
$$
$$\left\{\begin{aligned}
d(Y^v_t-\IE Y^v_t)=&\ \left[f_t(x^v_t-\IE x^v_t)+\tilde{g}_t(v_t-\IE v_t)\right]dt\\
                   &\ +h_td\tilde w_t,\\
Y^{v}_0-\IE Y^v_0=&\ 0,
\end{aligned}
\right.
$$
respectively. Let
$$
\mathbf x_0=\left(\begin{array}{c}
\xi-\mu_0\\
\mu_0
\end{array}
\right), \quad
\mathbf v_t=\left(\begin{array}{c}
v_t-\IE v_t\\
\IE v_t
\end{array}
\right),
$$
$$
\mathbf x^v_t=\left(\begin{array}{c}
x^v_t-\IE x^v_t\\
\IE x^v_t
\end{array}
\right), \quad
\mathbf y^v_t=\left(\begin{array}{c}
y^v_t-\IE y^v_t\\
\IE y^v_t
\end{array}
\right),
$$
$$
\mathbf z^v_t=\left(\begin{array}{c}
z^v_t-\IE z^v_t\\
\IE z^v_t
\end{array}
\right),
\quad
\mathbf{\tilde z}^v_t=\left(\begin{array}{c}
\tilde z^v_t-\IE \tilde z^v_t\\
\IE \tilde z^v_t
\end{array}
\right),
$$
$$
\mathbf Y^v_t=\left(\begin{array}{c}
Y^v_t-\IE Y^v_t\\
\IE Y^v_t
\end{array}
\right),
$$
and
$$
\mathbf a_t=\left(\begin{array}{cc}
a_t& 0\\
0&   a_t+\bar a_t
\end{array}
\right), \quad
\mathbf b_t=\left(\begin{array}{cc}
\tilde{b}_t& 0\\
0&   \tilde{b}_t
\end{array}
\right), \quad
\mathbf{\bar b}_t=\left(\begin{array}{c}
0\\
\bar b_t
\end{array}
\right), $$
$$
\mathbf c_t=\left(\begin{array}{c}
c_t\\
0
\end{array}
\right), \quad
\check\alpha_t=\left(\begin{array}{cc}
\alpha_t& 0\\
0&   \alpha_t+\bar \alpha_t
\end{array}
\right),
$$
$$
\check\beta_t=\left(\begin{array}{cc}
\beta_t& 0\\
0&   \beta_t+\bar\beta_t
\end{array}
\right),\quad
\check\gamma_t=\left(\begin{array}{cc}
\gamma_t& 0\\
0&   \gamma_t+\bar\gamma_t
\end{array}
\right),
$$
$$
\check{\tilde\gamma}_t=\left(\begin{array}{cc}
\tilde\gamma_t& 0\\
0&   \tilde\gamma_t+\bar{\tilde\gamma}_t
\end{array}
\right), \quad\check\psi_t=\left(\begin{array}{cc}
\tilde{\psi}_t& 0\\
0&   \tilde{\psi}_t
\end{array}
\right), \quad
\bar{\check \psi}_t=\left(\begin{array}{c}
0\\
\bar\psi_t
\end{array}
\right),$$
$$
\check c=\left(\begin{array}{cc}
1& 1\\
0& 0
\end{array}
\right), \quad
\check\rho=\left(\begin{array}{cc}
\rho & 0\\
0& \rho+\bar\rho
\end{array}
\right), \quad
\mathbf f_t=\left(\begin{array}{cc}
f_t& 0\\
0&   f_t+\bar f_t
\end{array}
\right),
$$
$$
\mathbf g_t=\left(\begin{array}{c}
0\\
\tilde{g}_t
\end{array}
\right), \quad
\mathbf h_t=\left(\begin{array}{c}
h_t\\
0
\end{array}
\right).
$$
Then
\begin{equation}
\label{cse}\left\{\begin{aligned}
d\mathbf x^v_t=&\ \left(\mathbf a_t\mathbf x^v_t+\mathbf b_t\mathbf v_t+\bar{\mathbf b}_t\right)dt+\mathbf c_tdw_t,\\
-d\mathbf y^v_t=&\ \left(\check\alpha_t\mathbf x^v_t+\check\beta_t\mathbf y^v_t+\check\gamma_t\mathbf z^v_t+\check{\tilde\gamma}_t\mathbf{\tilde z}^v_t+\check\psi_t\mathbf v_t\right.\\
                &\ \left.+\bar{\check\psi}_t\right)dt-\check c\mathbf z^v_tdw_t-\check c\mathbf{\tilde z}^v_td\tilde w_t,\\
\mathbf x^v_0=&\ \mathbf x_0, \quad\mathbf y_T=\check\rho\mathbf x^v_T,
\end{aligned}
\right.
\end{equation}
\begin{equation}
\label{coe}
\left\{\begin{aligned}
d\mathbf Y^v_t=&\ (\mathbf f_t\mathbf x^v_t+\mathbf g_t)dt+\mathbf h_td\tilde w_t,\\
\mathbf Y^v_0=&\ 0.
\end{aligned}
\right.
\end{equation}
On the other hand, let
$$
\mathbf A_t=\left(\begin{array}{cc}
A_t& 0\\
0&   A_t+\bar A_t
\end{array}
\right), \quad
\mathbf B_t=\left(\begin{array}{cc}
B_t& 0\\
0&   B_t
\end{array}
\right),$$
$$
\mathbf H=\left(\begin{array}{cc}
H& 0\\
0&   H+\bar H
\end{array}
\right), \quad
\mathbf M=\left(\begin{array}{cc}
M& 0\\
0&   M
\end{array}
\right).
$$
By simple calculations, we get
$$
\begin{aligned}
\IE\{\langle A_tx^v_t, x^v_t\rangle+\langle \bar A_t\IE x^v_t, \IE x^v_t\rangle\}=&\ \IE\langle\mathbf A_t\mathbf x^v_t, \mathbf x^v_t\rangle,\\
\IE\{\langle Hx^v_T, x^v_T\rangle+\langle \bar H\IE x^v_T, \IE x^v_T\rangle\}=&\ \IE\langle\mathbf H\mathbf x^v_T, \mathbf x^v_T\rangle,\\
\IE\langle My^v_0, y^v_0\rangle=\IE\langle\mathbf M\mathbf y^v_0, \mathbf y^v_0\rangle,\ \IE\langle B_tv_t, v_t\rangle=&\ \IE\langle\mathbf B_t\mathbf v_t, \mathbf v_t\rangle.
\end{aligned}
$$
Then the cost functional is rewritten as
\begin{equation}
\label{ccf}\begin{aligned}
J[v]=&\ \IE\left[\int_0^T(\langle\mathbf A_t\mathbf x^v_t, \mathbf x^v_t\rangle+\langle\mathbf B_t\mathbf v_t, \mathbf v_t\rangle)dt\right.\\
&\ \left.+\langle\mathbf H\mathbf x^v_T, \mathbf x^v_T\rangle+\langle\mathbf M\mathbf y^v_0, \mathbf y^v_0\rangle\right].
\end{aligned}
\end{equation}
Note that (\ref{ccf}) together with (\ref{cse}) and (\ref{coe}) forms a standard-looking LQ problem for FBSDE with noisy observation. However, the BSDE in (\ref{cse}) is not a standard form due to the irreversibility of $\check c$. Moreover, the control domain has to satisfy some extra constraint conditions according to the form of the control $\mathbf v$. This implies that Example (MFLQ) cannot be reduced to a standard LQ problem for FBSDE, hence it cannot be immediately solved by the standard LQ theory for FBSDE.

In the end of this section, we give a preliminary result, which shows that the desired optimality condition can be derived by minimizing $J[v]$ over $\mathcal U^0_{ad}$.

\textbf{Theorem 2.1.} $$\inf_{v'\in\mathscr U_{ad}}J[v']=\inf_{v\in\mathscr U^0_{ad}}J[v].$$

The proof can be found in Appendix.     \hfill$\Box$

\section{Optimal solution of Problem (MFC)}

For any $v$, $v_j\in\mathscr U_{ad}$, let $(x^v, y^v, z^v, \tilde{z}^v)$ and  $(x^{v_{j}}, y^{v_{j}}, z^{v_{j}},$ $\tilde{z}^{v_{j}})$ be the solutions of
(\ref{se}) corresponding to $v$ and $v_j$, $j=1, 2, \cdots$. For simplicity, we set
$$
\begin{aligned}
(\Theta^\lambda_t)=&\ (t, x^\lambda_t, \IE  x^\lambda_t, \lambda_t), \quad
(\Xi^\lambda_t)=(x^\lambda_t, \IE x^\lambda_t),\\
(\Pi_{t}^{\lambda})=& (x^{\lambda}_{t},  y^{\lambda}_{t},  z^{\lambda}_{t},  \tilde{z}^{\lambda}_{t}, \IE x^{\lambda}_{t}, \IE y^{\lambda}_{t}, \IE z^{\lambda}_{t}, \IE \tilde{z}^{\lambda}_{t})
\end{aligned}
$$ with $\lambda=v, u, v_j$, $j=1, 2, \cdots$.

\subsection{Optimality conditions}

According to Theorem 2.1 above, the optimality conditions can be derived by minimizing $J[v]$ over $\mathscr U^0_{ad}$ subject to (\ref{se}) and (\ref{oe}). We remind reader again that these results are different from the existing literature, mainly due to some new features of Problem (MFC). For example, the state is governed by a mean-field FBSDE, and is partially observed via a noisy process.

\textbf{Theorem 3.1.} \emph{Assume that $u$ is an optimal control for Problem (MFC).
Then the mean-field FBSDE
\begin{equation}\label{ae}
\left\{\begin{aligned}
dk_t=&\ \left(\beta^\top_t k_t+\bar{\beta}^\top_t\IE k_t\right)dt+\left(\gamma^\top_t k_t+\bar{\gamma}^\top_t\IE k_t\right)dw_t\\
     &\ +\left(\tilde{\gamma}^\top_t k_t+\bar{\tilde{\gamma}}^\top_t\IE k_t\right)d\tilde{w}_t,\\
-dp_t=&\ \left[a^\top_tp_t+l^\top_x(\Theta^u_t)+\IE\left(\bar a^\top_tp_t+l^\top_{\bar x}(\Theta^u_t)\right)-\alpha^\top_tk_t\right.\\
      &\ \left.-\bar{\alpha}^\top_t \IE k_t\right]dt-q_tdw_t-\tilde q_td\tilde w_t,\\
 k_0=&\ -\varphi^\top_{y}(y^{u}_0),\\
 p_T=&\ \phi^\top_x(\Xi^u_T)+\IE \phi^\top_{\bar x}(\Xi^u_T)-\rho^\top k_T-\bar{\rho}^\top\IE k_T
\end{aligned}
\right.
\end{equation}
has a unique solution $(k, p, q, \tilde q)$ in the space $\mathscr L^2_{\mathscr F}(0, T;$ $\IR^{m+n+n\times r+n\times\tilde r})$ such that
\begin{equation}
\label{mc}\IE\left[H_v(t, \Pi^u_t, u_t; k_t, p_t, q_t)(\nu-u_t)\big|\mathscr F^{Y^u}_t\right]\geq0
\end{equation}
for any $\nu\in U$, where the Hamiltonian function $H$
is defined by}
\begin{equation*}
\begin{aligned}
& H(t, x, y, z, \tilde{z}, \bar{x}, \bar{y}, \bar{z}, \bar{\tilde{z}}, v; k, p, q)\\
=& \langle a_tx+\bar a_t\bar x+b(t, v),p\rangle+\langle c_t, q\rangle+l(t, x, \bar x, v)\\
&-\left\langle \alpha_t x+ \bar{\alpha}_t \bar{x}+ \beta_t y+ \bar{\beta}_t \bar{y}+\gamma_t z+ \bar{\gamma}_t \bar{z}+ \tilde{\gamma}_t \tilde{z}\right.\\
&\left.+\bar{\tilde{\gamma}}_t \bar{\tilde{z}}+\psi(t, v), k\right\rangle.
\end{aligned}
\end{equation*}

\textbf{Proof.} If $u$ is an optimal control for Problem (MFC), Theorem 2.1 implies that $J[u]=\inf_{v\in\mathscr U^0_{ad}}J[v].$ For any
$v\in\mathscr U^0_{ad}$, let $(x^{u+\varepsilon v}, y^{u+\varepsilon v}, z^{u+\varepsilon v}, \tilde{z}^{u+\varepsilon v})$ be the solution of (\ref{se}) corresponding to $u+\varepsilon v$, $0<\varepsilon<1$. Introduce a variational equation
$$
\left\{\begin{aligned}
\dot x^u_{1, t}=&\ a_tx^u_{1, t}+\bar a_t\IE x^u_{1, t}+b_v(t,u_t)v_t, \\
-dy^u_{1, t}=&\ \left(\alpha_tx^u_{1, t}+\bar\alpha_t\IE x^u_{1, t}+\beta_ty^u_{1, t}+\bar\beta_t\IE y^u_{1, t}\right.\\
&\ +\gamma_tz^u_{1, t}+\bar\gamma_t\IE z^u_{1, t}+\tilde\gamma_t\tilde z^u_{1, t}+\bar{\tilde\gamma}_t\IE\tilde z^u_{1, t}\\
&\ \left.+\psi_{v}(t, u_t)v_{t}\right)dt-z^u_{1, t}dw_t-\tilde z^u_{1, t}d\tilde w_t,\\
x^u_{1, 0}=&\ 0, \quad y^u_{1, T}=\rho x^u_{1, T}+\bar\rho\ \IE x^u_{1, T},
\end{aligned}
\right.
$$
which admits a unique solution $(x^u_1, y^u_1, z^u_1, \tilde{z}^u_1) \in\mathscr L^2_{\mathscr F}(0, T; \IR^{n+m+m\times r+m\times \tilde{r}})$. It follows from Taylor's expansion, H\"{o}lder's inequality and the techniques applied in Lemma A.1 (see Appendix) that
$$\lim_{\varepsilon\rightarrow0}\IE\sup_{0\leq t\leq T}\left|\frac{x^{u+\varepsilon v}_t-x^u_t}{\varepsilon}-x^u_{1,t}\right|^2=0,$$
and then,
$$\lim_{\varepsilon\rightarrow0}\IE\sup_{0\leq t\leq T}\left|\frac{y^{u+\varepsilon v}_t-y^u_t}{\varepsilon}-y^u_{1,t}\right|^2=0.$$
Combining the limits with the optimality of $u$ using the first variation of $J[v]$, we have
\begin{eqnarray}
\label{vin}
0&\leq&\lim_{\varepsilon\rightarrow0}\frac{J[u+\varepsilon v]-J[u]}{\varepsilon}\nonumber\\
&=&\IE\int_0^T\big(l^\top_x(\Theta^u_t)x^u_{1,t}+l^\top_{\bar x}(\Theta^u_t)\IE x^u_{1,t}+l^\top_v(\Theta^u_t)v_t\big)dt\nonumber\\
&&+\IE\big(\phi^\top_x(\Xi^u_T)x^u_{1, T}+\phi^\top_{\bar
x}(\Xi^u_T)\IE x^u_{1, T}+\varphi^\top_y(y^u_0)y^u_{1,0}\big).\nonumber\\
\end{eqnarray}
On the other hand, once $(x^u, y^u, z^u, \tilde{z}^u)$ is determined by (\ref{se}), there is
a unique solution $(k, p, q, \tilde q)$, in the space $\mathscr L^2_{\mathscr F}(0, T; \IR^{m+n+n\times r+n\times\tilde r})$, to
(\ref{ae}). Using It\^o's formula to $\langle x^u_{1}, p\rangle+\langle y^u_{1}, k\rangle$ and
inserting it into (\ref{vin}), we get
$$\IE\int_0^T\left(p^\top_t b_v(t,u_t)-k^\top_t \psi_v(t,u_t)+l^\top_v(\Theta^u_t)\right)v_tdt\geq0.$$
Due to $u\in\mathscr U^0_{ad}$ and the arbitrariness of $v_t$, we
deduce
\begin{align*}
 & \IE\left[\left(p^\top_t b_v(t,u_t)-k^\top_t \psi_v(t,u_t)+l^\top_v(\Theta^u_t)\right)(\nu-u_t)\big|\mathscr F^{Y^0}_t\right]\\
=&\, \IE\left[H_v(t, \Pi^u_t, u_t; k_t, p_t, q_t)(\nu-u_t)\big|\mathscr F^{Y^0}_t\right]\geq0,
\end{align*}
for any $\nu\in U$. Since $u\in\mathscr U_{ad}$, it follows from Proposition 2.1 that $\mathscr F^{Y^u}_t=\mathscr
F^{Y^0}_t$. Then the result is derived. The proof is complete.
\hfill$\Box$

\textbf{Theorem 3.2.} \emph{Assume that for any $(t, x, y, z, \tilde{z}, \bar x, \bar{y}, \bar{z}, \bar{\tilde{z} },$ $v)\in
[0, T]\times\IR^{n+m+m\times r+m\times \tilde{r}+n+m+m\times r+m\times \tilde{r}}\times U$, $(x, \bar x,
v)\mapsto l(t, x, \bar x, v)$, $(x, \bar x)\mapsto\phi(x, \bar
x)$ and $y\mapsto \varphi(y)$ are convex. Assume that $b(t, v_t)=\tilde{b}_tv_t+\bar b_t, \psi(t, v_t)=\tilde{\psi}_tv_t+\bar\psi_t,$ where
$\tilde{b}\in\mathscr L^\infty(0, T; \IR^{n\times k})$,
$\bar b\in\mathscr L^\infty(0, T; \IR^n)$, $\tilde{\psi}\in\mathscr L^\infty(0, T; \IR^{m\times k})$,
$\bar \psi\in\mathscr L^\infty(0, T; \IR^m)$. Let $u\in\mathscr U_{ad}$ and
\begin{equation}
\label{mc2}
\begin{aligned}
&\ \IE\left[H(t, \Pi^u_t, u_t; k_t, p_t, q_t)\big|\mathscr F^{Y^u}_t\right]\\
=&\ \inf_{v\in U}\IE\left[H(t, \Pi^u_t, v; k_t, p_t, q_t)\big|\mathscr F^{Y^u}_t\right],
\end{aligned}
\end{equation}
where $(x^u, y^u, z^u, \tilde{z}^u)\in\mathscr L^2_{\mathscr F}(0, T; \IR^{n+m+m\times r+m\times \tilde{r}})$ is the solution to \eqref{se} under the admissible control $u$, and $(k, p, q, \tilde q)\in\mathscr L^2_{\mathscr F}(0, T; \IR^{m+n+n\times r+n\times\tilde r})$ is the solution to \eqref{ae}.
Then $u$ is an optimal control of Problem (MFC).}

\textbf{Proof.} For any $v\in\mathscr U_{ad}$, we write
\begin{equation}
\label{tv}J[v]-J[u]=I_1+I_2+I_3
\end{equation}
with
\begin{align*}
I_1=&\ \IE\int_0^T\left(l(\Theta^v_t)-l(\Theta^u_t)\right)dt, I_2=\IE\left[\phi(\Xi^v_T)-\phi(\Xi^u_T)\right], \\
I_3=&\ \IE\left[\varphi(y^v_0)-\varphi(y^u_0)\right].
\end{align*}
By virtue of the convexity of $\phi$ and applying It\^o's formula to
$\langle p, x^v-x^u\rangle$, we have
\begin{eqnarray}\label{i2}
I_2&\geq&\IE\left\langle p_T, x^v_T-x^u_T\right\rangle+\IE\left\langle \rho^{\top}k_T+\IE(\bar{\rho}^{\top}k_T), x^v_T-x^u_T\right\rangle\nonumber\\
&=&-\IE\int_0^T\left\langle l^{\top}_x(\Theta^u_t)+\IE\, l^{\top}_{\bar x}(\Theta^u_t)-\alpha^{\top}_tk_t-\IE(\bar{\alpha}^{\top}_tk_t),\right.\nonumber\\
&& \left.\quad\quad\quad\quad x^v_t-x^u_t\right\rangle dt+\IE\int_0^T\left\langle p_t, \tilde{b}_{t}(v_t-u_t)\right\rangle dt\nonumber\\
&& +\IE\left\langle \rho^{\top}k_T+\IE(\bar{\rho}^{\top}k_T), x^v_T-x^u_T\right\rangle.
\end{eqnarray}
Similarly, applying It\^o's formula to
$\langle k, y^v-y^u\rangle$ and the convexity of $\varphi$, we derive
\begin{eqnarray}\label{i3}
  I_3 &\geq& -\IE\langle k_T,\; \rho(x^v_T-x^u_T)+\bar{\rho}\IE(x^v_T-x^u_T)\rangle\nonumber\\
  && -\IE\int_0^T\left\langle k_t, \alpha_t(x^v_t-x^u_t)+\bar{\alpha}_t\IE(x^v_t-x^u_t)\right.\nonumber\\
  && \left.\quad\quad\quad\quad+\tilde{\psi}_{t}(v_t-u_t)\right\rangle dt.
\end{eqnarray}
It is easy to see from (\ref{tv}), (\ref{i2}), (\ref{i3}) and the convexity of $l$ that
$$
\begin{aligned}
& J[v]-J[u]\\
\geq&\ \IE\int_0^T H_v\big(t, \Pi^u_t, u_t; k_t, p_t, q_t\big)\big(v_t-u_t\big)dt\\
=&\ \IE\int_0^T\IE\left[ H_v\big(t, \Pi^u_t, u_t; k_t, p_t, q_t\big)\big(v_t-u_t\big)\big|\mathscr F^{Y^0}_t\right]dt.
\end{aligned}
$$
Further, using Theorem  2.1 and (\ref{mc2}), we get
$$
\begin{aligned}
&\IE\left[ H_v\big(t, \Pi^u_t, u_t; k_t, p_t, q_t\big)\big(v_t-u_t\big)\big|\mathscr F^{Y^u}_t\right]\\
=& \frac{\partial}{\partial v}\IE\left[H\big(t, \Pi^u_t, u_t; k_t, p_t, q_t\big)\big|\mathscr F^{Y^u}_t\right](v_t-u_t)\geq 0. \\
\end{aligned}
$$
Then this implies the desired result.  \hfill$\Box$

\subsection{Optimal filters}

The minimum condition (\ref{mc}) shows that we need to analyze the
optimal filters of (\ref{se}) and (\ref{ae}) depending on
$\mathscr F^{Y^v}_t$ in order to compute $u$. For this, for any $v\in\mathscr U_{ad}$ we denote
by $\hat\iota_t=\IE\left[\iota_t\big|\mathscr F^{Y^v}_t\right]$ with $\iota_t=x^0_t, x^v_t,  y^v_t, z^v_t, \tilde{z}^v_t, y_t^{v}(x^v_t)^\top$
and $\hat\kappa_t=\IE\left[\kappa_t\big|\mathscr F^{Y^u}_t\right]$ with
$\kappa_t=k_t, p_t, p_t(x^u_t)^\top, k_t(x^u_t)^\top, \tilde q_t$ the optimal filters of
$\iota_t$ and $\kappa_t$, respectively. Moreover, we denote by
$\Sigma_t=\IE\left[(x^v_t-\hat x^v_t)(x^v_t-\hat x^v_t)^{\top}\right]$ the mean square error of $\hat x^v_t$, $v\in\mathscr U_{ad}$.

Using Theorem 12.7 in Liptser and Shiryayev \cite{ls} and Theorem 3.1 in Wang et al. \cite{wzz}, we derive (\ref{sfe}) and (\ref{afe}), respectively.

\textbf{Theorem 3.3.} \emph{For any $v\in\mathscr U_{ad}$, the optimal
filters $(\hat{x}^v_t, \hat{y}^v_t)$ and $(\hat{k}^v_t,\hat{ p}^v_t)$ of the solutions $(x^v_t, y^v_t)$ and $(k_t, p_t)$ to (\ref{se}) and
(\ref{ae}) with respect to $\mathscr F^{Y^v}_t$ and $\mathscr
F^{Y^u}_t$ satisfy
\begin{equation}
\label{sfe}\left\{\begin{aligned}
d\hat x^v_t=&\ \left(a_t\hat x^v_t+\bar a_t\IE x^v_t+b(t, v_t)\right)dt\\
            &\ +\Sigma_tf^\top_t(h_t^{-1})^{\top}d\bar w_t,\\
-d\hat{y}^v_t=&\ \left(\alpha_t\hat{x}^v_t+\bar\alpha_t\IE x^v_t+\beta_t\hat{y}^v_t+\bar\beta_t\IE y^v_t+\gamma_t\hat{z}^v_t\right.\\
&\ \left.+\bar\gamma_t\IE z^v_t+\tilde\gamma_t\hat{\tilde z}^v_t+\bar{\tilde\gamma}_t\IE\tilde z^v_t+\psi(t, v_t)\right)dt\\
&\ -\widehat{\tilde{Z}}_t^vd\bar w_t,\\
\hat x^v_0=&\ \mu_0, \quad \hat{y}^v_T=\rho\, \hat{x}^v_T+\bar\rho\, \IE x^v_T,
\end{aligned}
\right.
\end{equation}
and
\begin{equation}
\label{afe}\left\{\begin{aligned}
d\hat{k}_t=&\ \left(\beta^\top_t \hat{k}_t+\bar{\beta}^\top_t\IE k_t\right)dt+\left[\tilde{\gamma}^\top_t \hat{k}_t+\bar{\tilde{\gamma}}^\top_t \IE k_t\right.\\
&\ \left.+\left(\widehat{k_t (x^u_t)^\top}-\hat{k}_{t} (\hat{x}^u_t)^{\top}\right)f^\top_t(h^{-1}_t)^\top\right]d\bar{w}_t,\\
-d\hat p_t=&\ \left\{a^\top_t\hat p_t+\IE\left[l_x^\top(\Theta^u_t)\big|\mathscr F^{Y^u}_t\right]-\alpha^\top_t\hat{k}_t-\bar{\alpha}^\top_t\IE k_t\right.\\
  &\ \left.+\IE\left(\bar a^\top_tp_t+l_{\bar x}^\top(\Theta^u_t)\right)\right\}dt-\widehat{\tilde{Q}}_td\bar w_t,\\
\hat{k}_0=&\ -\varphi^\top_{y}(\hat{y}^{u}_0),\\
\hat p_T=&\ \IE\left[\phi_x^{\top}(\Xi^u_T)\big|\mathscr F^{Y^u}_T\right]+\IE\phi_{\bar x}^{\top}(\Xi^u_T)\\
         &\ -\rho^\top \hat{k}_T-\bar{\rho}^\top\IE k_T,
\end{aligned}
\right.
\end{equation}
respectively, where $\Sigma$ is the unique solution of
\begin{equation}
\label{sigma}\left\{\begin{aligned}
&\ \dot\Sigma_t-a_t\Sigma_t-\Sigma_ta^\top_t+\Sigma_tf^\top_t(h^{-1}_t)^{\top}h^{-1}_tf_t\Sigma_t\\
&\ -c_tc^\top_t=0,\\
&\ \Sigma_0=\sigma_0,
\end{aligned}
\right.
\end{equation}
\begin{eqnarray}
\label{ip}\bar w_t&=&\int_0^th^{-1}_s\left[dY^0_s-\left(f_s\hat
x^0_s+\bar f_s\IE x^0_s\right)ds\right]
\end{eqnarray}
is a standard Brownian motion with value in $\IR^{\tilde r}$, and}
$$\widehat{\tilde{Z}}_t=\hat{\tilde{z}}_t+\left(\widehat{y_t^{v} (x^v_t)^\top}-\hat{y}_{t}^{v} (\hat{x}^v_t)^{\top}\right)f^\top_t(h^{-1}_t)^\top,$$
$$\widehat{\tilde{Q}}_t=\hat{\tilde q}_t+\left(\widehat{p_t (x^u_t)^\top}-\hat{p}_{t} (\hat{x}^u_t)^{\top}\right)f^\top_t(h^{-1}_t)^\top.$$

We emphasize that (\ref{sfe}) and (\ref{afe}) are two forward-backward optimal filters. It shows that the difference between Theorem 3.3 and the classical filtering literature, say, Bensoussan \cite{ba1,ba2}, Liptser and Shiryayev \cite{ls}.

\section{An LQ case of Problem (MFC)}

We still adopt the notations and the assumptions introduced in Sections 2 and 3 unless noted otherwise.

\textbf{Problem (MFLQ).} Minimize
\begin{align}\label{lcf}
J[v]=&\ \frac{1}{2}\IE\left\{\int_0^T\left[\langle A_tx^v_t, x^v_t\rangle+\langle \bar A_t\IE x^v_t, \IE x^v_t\rangle\right.\right.\nonumber\\
&\ +\langle B_tv_t, v_t\rangle+2\langle D_tx^v_t,v_t\rangle+2\langle\bar D_t\IE x^v_t, v_t\rangle\nonumber\\
&\ \left.+2\langle \tilde{F}_t,x^v_t\rangle+2\langle\tilde{\bar{F}}_{t}, \IE x^v_t\rangle+2\langle \tilde{G}_t, v_t\rangle\right]dt\nonumber\\
&\ +\langle Hx^v_T, x^v_T\rangle+\langle\bar H\IE x^v_T, \IE x^v_T\rangle+2\langle \tilde{L}, x^v_T\rangle\nonumber\\
&\ \left.+2\langle\tilde{\bar L}, \IE x^v_T\rangle+\langle My^v_0, y^v_0\rangle+2\langle N, y^v_0\rangle\right\}
\end{align}
over $\mathscr U_{ad}$ with the control domain $U=\IR^k$, subject to the state equation
\begin{equation}\label{se3}
\left\{\begin{aligned}
dx^v_t=&\ \left(a_tx^v_t+\bar a_t\IE x^v_t+b_tv_t+\bar b_t\right)dt+c_tdw_t,\\
-dy^v_t=&\ \left(\alpha_tx^v_t+\bar\alpha_t\IE x^v_t+\beta_ty^v_t+\bar\beta_t\IE y^v_t+\gamma_tz^v_t\right.\\
&\ \left.+\bar\gamma_t\IE z^v_t+\tilde\gamma_t\tilde z^v_t+\bar{\tilde\gamma}_t\IE\tilde z^v_t+\psi_{t}v_{t}+\bar{\psi_{t}}\right)dt\\
&\ -z^v_tdw_t-\tilde z^v_td\tilde w_t,\\
x^v_0=&\ \xi, \quad y^v_T=\rho x^v_T+\bar\rho\, \IE x^v_T,
\end{aligned}
\right.
\end{equation}
and the observation equation
\begin{equation}\label{oe3}
\left\{\begin{aligned}
dY^v_t=&\ \left(f_tx^v_t+\bar f_t\IE x^v_t+g_t\right)dt+h_td\tilde w_t,\\
Y^{v}_0=&\ 0,
\end{aligned}
\right.
\end{equation}
where $A, \bar A\in\mathscr
L^\infty(0, T; S^n)$, $B\in\mathscr L^\infty(0, T; S^k)$, $H, \bar
H\in S^n, $ $A, H, M \geq0$, $B>0$, $A+\bar A, H+\bar H \geq0$, $D, \bar
D\in\mathscr L^\infty(0, T; \IR^{k\times n})$, $\tilde F, \tilde{\bar F}\in\mathscr
L^\infty(0, T; \IR^n)$, $\tilde G\in\mathscr L^\infty(0, T; \IR^k)$, $L,
\bar L\in\IR^n, N\in\IR^m$, $b\in\mathscr L^\infty(0, T; \IR^{n\times k})$,
$\bar b\in\mathscr L^\infty(0, T; \IR^n)$, $\psi\in\mathscr L^\infty(0, T; \IR^{m\times k})$,
$\bar \psi\in\mathscr L^\infty(0, T; \IR^m)$, and $g\in\mathscr L^\infty(0, T; \IR^{\tilde r})$.

Note that we do not assume the positive semidefiniteness of $\bar A$ and $\bar H$ in Problem (MFLQ). Then (\ref{lcf}) covers the performance functional of Problem (AL) as a special case.

\textbf{Proposition 4.1.} \emph{If $u$ is an optimal control for Problem (MFLQ), then
$$\begin{aligned}
u_t=&\ -B^{-1}_t\left(b^\top_t\IE\left[p_t|\mathscr F^{Y^u}_t\right]-\psi^\top_t\IE\left[k_t|\mathscr F^{Y^u}_t\right]\right.\\
    &\ \left.+D_t\IE\left[x^u_t|\mathscr F^{Y^u}_t\right]+\bar D_t\IE x^u_t+\tilde G_t\right),
\end{aligned}
$$
where $(k, p)$ is the solution of the adjoint equation
\begin{equation}\label{lae1}
\left\{\begin{aligned}
dk_t=&\ \left(\beta^\top_t k_t+\bar{\beta}^\top_t\IE k_t\right)dt+\left(\gamma^\top_t k_t+\bar{\gamma}^\top_t\IE k_t\right)dw_t\\
     &\ +\left(\tilde{\gamma}^\top_t k_t+\bar{\tilde{\gamma}}^\top_t\IE k_t\right)d\tilde{w}_t,\\
-dp_t=&\ \left[a^\top_tp_t-\alpha_{t}^{\top}k_t+A_tx^u_t+D^\top_tu_t+F_t\right.\\
      &\ \left.+\IE\left(\bar a^\top_tp_t-\bar{\alpha}_{t}^{\top}k_t+\bar A_tx^u_t+\bar D^\top_tu_t+\bar F_t\right)\right]dt\\
      &\ -q_tdw_t-\tilde q_td\tilde w_t,\\
k_0=&\ -My^{u}_0-N, \\
p_T=&\ Hx^u_T+\bar H\IE x^u_T+\tilde L+\tilde{\bar L}-\rho^\top k_T-\bar{\rho}^\top\IE k_T,
\end{aligned}
\right.
\end{equation}
together with the state equation
\begin{equation}
\left\{\begin{aligned}\label{olse1}
dx^u_t=&\ \left(a_tx^u_t+\bar a_t\IE x^u_t+b_tu_t+\bar b_t\right)dt+c_tdw_t,\\
-dy^u_t=&\ \left(\alpha_tx^u_t+\bar\alpha_t\IE x^u_t+\beta_ty^u_t+\bar\beta_t\IE y^u_t+\gamma_tz^u_t\right.\\
        &\ \left.+\bar\gamma_t\IE z^u_t+\tilde\gamma_t\tilde z^u_t+\bar{\tilde\gamma}_t\IE\tilde z^u_t+\psi_{t}u_{t}+\bar{\psi_{t}}\right)dt\\
        &\ -z^u_tdw_t-\tilde z^u_td\tilde w_t,\\
x^u_0=&\ \xi, \quad y^u_T=\rho x^u_T+\bar\rho\ \IE x^u_T.
\end{aligned}
\right.
\end{equation}
}

\textbf{Proof.} With the above data, the Hamiltonian function is
$$
\begin{aligned}
 &\ H(t, x, y, z, \tilde{z}, \bar{x}, \bar{y}, \bar{z}, \bar{\tilde{z}}, v; k, p, q)\\
=&\ \langle a_tx+\bar a_t\bar x+b_tv+\bar b_t, p\rangle+\langle c_t, q\rangle-\left\langle \alpha_t x+ \bar{\alpha}_t \bar{x}+ \beta_t y\right.\\
&\ \left.+\bar{\beta}_t \bar{y}+\gamma_t z+ \bar{\gamma}_t \bar{z}+ \tilde{\gamma}_t \tilde{z}+ \bar{\tilde{\gamma}}_t \bar{\tilde{z}}+\psi_{t}v+\bar{\psi}_{t}, k\rangle\right.\\
&\ +\frac{1}{2}\left[\langle A_tx, x\rangle+\langle \bar A_t\bar x, \bar x\rangle+\langle B_tv, v\rangle +2\langle D_tx,v\rangle\right.\\
&\ \left.+2\langle\bar D_t\bar x, v\rangle+2\langle \tilde F_t, x\rangle+2\langle\tilde{\bar F}_t, \bar x\rangle+2\langle \tilde G_t, v\rangle\right],
\end{aligned}
$$
where $(k, p, q, \tilde{q})$ is determined by (\ref{lae1}) together with (\ref{olse1}). If $u(\cdot)$ is optimal, it follows from Theorem 3.1 that
$$
\begin{aligned}
u_t=&\ -B^{-1}_t\left(b^\top_t\IE\left[p_t|\mathscr F^{Y^u}_t\right]-\psi^\top_t\IE\left[k_t|\mathscr F^{Y^u}_t\right]\right.\\
&\ \left.+D_t\IE\left[x^u_t|\mathscr F^{Y^u}_t\right]+\bar D_t\IE x^u_t+\tilde G_t\right),
\end{aligned}
$$
where $(x^u, k, p)$ is the solution of (\ref{olse1}) with (\ref{lae1}). Then the proof is complete. \hfill$\Box$

Note that one more explicit optimal control $u$ strong depends on a certain special structure of the state equation and the cost functional. Next, let us consider a particular case of Problem (MFLQ), i.e., let $M=0$ and $\bar\beta_t=\gamma_t=\bar\gamma_t=\tilde\gamma_t=\bar{\tilde\gamma}_t=0$ in (\ref{lcf}) and (\ref{se3}), respectively. By Theorems 3.1, 3.2 and 3.3, an optimal feedback control is explicitly obtained. The procedure of how to solve is decomposed into five steps below. Note that such an optimal control will play a role in Problem (AL). Please refer to Section 5 below for more details.

\emph{Step 1:} A reduced LQ problem.

Integrating and taking expectations on both sides of the BSDE in \eqref{se3}, we have
$$
\begin{aligned}
\IE y^{v}_{t}=&\ \IE \left[\chi_{t}^{T}\rho x_{T}^{v}+\chi_{t}^{T}\bar{\rho}\IE x_{T}^{v}\right.\\
              &\ \left.+\int_t^T\chi_{t}^{s}\left(\alpha_{s}x^{v}_{s}+\bar{\alpha}_{s}\IE x_{s}^{v}+\psi_{s}v_s+\bar{\psi}_{s}\right)ds\right]
\end{aligned}
$$
with $$\chi_{t}^{s}=e^{\int_{t}^{s}(\beta_\tau+\bar{\beta}_\tau)d\tau}, \quad t\leq s\leq T.$$
Plugging the equality into (\ref{lcf}), we derive an LQ problem for mean-field SDE as follows.

\textbf{Problem (MFLQ)$'$.} Find a $v\in\mathscr U_{ad}$ to minimize
$$\begin{aligned}
J[v]=&\ \frac{1}{2}\IE\left\{\int_0^T\left[\langle A_tx^v_t, x^v_t\rangle+\langle \bar A_t\IE x^v_t, \IE x^v_t\rangle\right.\right.\nonumber\\
&\ +\langle B_tv_t, v_t\rangle+2\langle D_tx^v_t,v_t\rangle+2\langle\bar D_t\IE x^v_t, v_t\rangle\\
&\ \left.+2\langle F_t,x^v_t\rangle+2\langle\bar F_t, \IE x^v_t\rangle+2\langle G_t, v_t\rangle\right]dt\nonumber\\
&\ +\langle Hx^v_T, x^v_T\rangle+\langle\bar H\IE x^v_T, \IE x^v_T\rangle+2\langle L, x^v_T\rangle\\
&\ \left.+2\langle\bar L, \IE x^v_T\rangle\right\}+J_{0}
\end{aligned}
$$
subject to
$$
\left\{\begin{aligned}
dx^v_t=&\ \left(a_tx^v_t+\bar a_t\IE x^v_t+b_tv_t+\bar b_t\right)dt+c_tdw_t,\\
x^v_0=&\ \xi
\end{aligned}
\right.
$$
and \eqref{oe3} with
\begin{align}
F&=\tilde{F}+(\chi_{0}^{t}\alpha_t)^{\top} N, \hspace{1.8cm}   \bar{F}=\tilde{\bar{F}}+(\chi_{0}^{t}\bar{\alpha}_t)^{\top} N,\nonumber\\
L&=\tilde{L}+(\chi_{0}^{T}\rho)^{\top} N, \hspace{2cm}  \bar{L}=\tilde{\bar{L}}+(\chi_{0}^{T}\bar{\rho})^{\top} N,\nonumber\\
G&=\tilde{G}+(\chi_{0}^{t}\psi_t)^{\top} N, \hspace{1.75cm} J_{0}=\int_0^TN^{\top}\chi_{0}^{t}\bar{\psi}_{t}dt.\nonumber
\end{align}

\emph{Step 2:} Candidate optimal control.

The Hamiltonian function is
$$
\begin{aligned}
H(t, x, \bar x, v; p, q)=&\ \langle a_tx+\bar a_t\bar x+b_tv+\bar b_t, p\rangle+\langle c_t, q\rangle\\
&\ +\frac{1}{2}[\langle A_tx, x\rangle+\langle \bar A_t\bar x, \bar x\rangle+\langle B_tv, v\rangle\\
&\ +2\langle D_tx,v\rangle+2\langle\bar D_t\bar x, v\rangle+2\langle F_t, x\rangle\\
&\ +2\langle\bar F_t, \bar x\rangle+2\langle G_t,
v\rangle],
\end{aligned}
$$
where $(p, q)$ is determined by the Hamiltonian system
\begin{equation}
\label{lhs}\left\{\begin{aligned}
dx^u_t=&\ \left(a_tx^u_t+\bar a_t\IE x^u_t+b_tu_t+\bar b_t\right)dt+c_tdw_t,\\
-dp_t=&\ \left[a^\top_tp_t+A_tx^u_t+D^\top_tu_t+F_t\right.\\
&\ \left.+\IE\left(\bar a^\top_tp_t+\bar A_tx^u_t+\bar D^\top_tu_t+\bar F_t\right)\right]dt\\
&\ -q_tdw_t-\tilde q_td\tilde w_t,\\
x^u_0=&\ \xi, \quad p_T=Hx^u_T+\bar H\IE x^u_T+L+\bar L.
\end{aligned}
\right.
\end{equation}
If $u$ is optimal, then it follows from Theorem 3.1 or Proposition 4.1 that
\begin{align}
\label{loc}
u_t=&\ -B^{-1}_t\left(b^\top_t\IE\left[p_t|\mathscr F^{Y^u}_t\right]+D_t\IE\left[x^u_t|\mathscr F^{Y^u}_t\right]\right.\nonumber\\
    &\ \left.+\bar D_t\IE x^u_t+G_t\right)\nonumber\\
=&\ -B^{-1}_t\left(b^\top_t\hat p_t+D_t\hat x^u_t+\bar D_t\IE x^u_t+G_t\right).
\end{align}

\emph{Step 3:} Feedback representation of (\ref{loc}).

Inserting (\ref{loc}) into (\ref{lhs}) and taking expectations, we
get a fully coupled forward-backward ordinary differential equation (ODE, in short)
\begin{equation}
\label{lhsm}\left\{\begin{aligned}
\dot{\IE x^u_t}=&\ \left[a_t+\bar a_t-b_tB^{-1}_t\left(D_t+\bar D_t\right)\right]\IE x^u_t\\
                &\ -b_tB^{-1}_tb^\top_t\IE p_t-b_tB^{-1}_tG_t+\bar b_t,\\
\dot{\IE p_t}=&\ -\left[A_t+\bar A_t-\left(D_t+\bar D_t\right)^\top B^{-1}_t(D_t+\bar D_t)\right]\IE x^u_t\\
              &\ -\left[(a_t+\bar a_t)^\top-\left(D_t+\bar D_t\right)^\top B^{-1}_tb^\top_t\right]\IE p_t\\
              &\ +(D_t+\bar D_t)^\top B^{-1}_tG_t-F_t-\bar F_t,\\
\IE x^u_0=&\ \mu_0, \quad \IE p_T=(H+\bar H)\IE x^u_T+L+\bar L.
\end{aligned}
\right.
\end{equation}
According to Theorem 2.6 in Peng and Wu \cite{pw}, (\ref{lhsm}) has a unique solution $(\IE x^u, \IE p)$ based on the assumption

\noindent \textbf{(A1).} There is a constant $C\geq0$ such that
$$CI_{n\times n}-\left[A_t+\bar A_t-\left(D_t+\bar D_t\right)^\top
B^{-1}_t\left(D_t+\bar D_t\right)\right]\leq0.$$
Hereinafter, $I_{n\times n}$ stands for an $n\times n$ unit matrix.

Noticing the terminal condition of (\ref{lhsm}), we set
\begin{equation}
\label{pxmr}\IE p_t=\Phi_t\IE x^u_t+\Psi_t
\end{equation}
for two deterministic and differentiable functions $\Phi$ and
$\Psi$ such that $\Phi_T=H+\bar H$ and $\Psi_T=L+\bar L$.
Applying the chain rule for computing the derivative of (\ref{pxmr}), we have
$$
\begin{aligned}
\dot{\IE p_t}=&\ \dot\Phi_t\IE x^u_t+\Phi_t\dot{\IE x^u_t}+\dot\Psi_t\\
=&\ \left\{\dot\Phi_t+\Phi_t\left[a_t+\bar a_t-b_tB^{-1}_t\left(D_t+\bar D_t\right)\right]\right.\\
&\ \left.-\Phi_tb_tB^{-1}_tb^\top_t\Phi_t\right\}\IE x^u_t+\dot\Psi_t-\Phi_tb_tB^{-1}_tb^\top_t\Psi_t\\
&\ +\Phi_t\left(\bar b_t-b_tB^{-1}_tG_t\right).
\end{aligned}
$$
Comparing it with the second equation in (\ref{lhsm}), we deduce a Ricatti equation
\begin{equation}
\label{alpha}\left\{\begin{aligned}
&\ \dot\Phi_t+\Phi_t\left[a_t+\bar a_t-b_tB^{-1}_t\left(D_t+\bar D_t\right)\right]\\
&\ +\left[(a_t+\bar a_t)^\top-\left(D_t+\bar D_t\right)^\top B^{-1}_tb^\top_t\right]\Phi_t\\
&\ -\Phi_tb_tB^{-1}_tb^\top_t\Phi_t+A_t+\bar A_t\\
&\ -\left(D_t+\bar D_t\right)^\top B^{-1}_t\left(D_t+\bar D_t\right)=0,\\
&\ \Phi_T=H+\bar H
\end{aligned}
\right.
\end{equation}
and an ODE
\begin{equation}
\label{beta}\left\{\begin{aligned}
&\ \dot\Psi_t+\left[\left(a_t+\bar a_t\right)^\top-\left(D_t+\bar D_t\right)^\top B^{-1}_tb^\top_t\right.\\
&\ \left.-\Phi_tb_tB^{-1}_tb^\top_t\right]\Psi_t+\Phi_t\left(\bar b_t-b_tB^{-1}_tG_t\right)\\
&\ -\left(D_t+\bar D_t\right)^\top B^{-1}_tG_t+F_t+\bar F_t=0,\\
&\ \Psi_T=L+\bar L.
\end{aligned}
\right.
\end{equation}
Clearly, (\ref{alpha}) admits a unique solution, and thus, (\ref{beta}) also has a unique solution.
Plugging (\ref{pxmr}) into the first equation of (\ref{lhsm}), we derive
\begin{equation}
\label{lhxm}\left\{\begin{aligned}
\dot{\IE x^u_t}=&\ \left[a_t+\bar a_t-b_tB^{-1}_t\left(D_t+\bar D_t\right)\right.\\
                &\ \left.-b_tB^{-1}_tb^\top_t\Phi_t\right]\IE x^u_t-b_tB^{-1}_tb^\top_t\Psi_t\\
                &\ -b_tB^{-1}_tG_t+\bar b_t,\\
\IE x^u_0=&\ \mu_0,
\end{aligned}
\right.
\end{equation}
which can be explicitly computed.

Using Theorem 3.3 to (\ref{lhs}) with (\ref{loc}), we get the
optimal filtering equation
\begin{equation}
\label{lhsf}\left\{\begin{aligned}
d\hat x^u_t=&\ \left[\left(a_t-b_tB^{-1}_tD_t\right)\hat x^u_t-b_tB^{-1}_tb^\top_t\hat p_t+\theta_{1,t}\right]dt\\
            &\ +\Sigma_tf_t^{\top}(h^{-1}_t)^{\top}d\bar w_t,\\
-d\hat p_t=&\ \left[\left(A_t-D^\top_tB^{-1}_tD_t\right)\hat x^u_t+\left(a^\top_t-D^\top_tB^{-1}_tb^\top_t\right)\hat p_t\right.\\
           &\ \left.+\theta_{2,t}\right]dt-\widehat{ \tilde{Q}}_td\bar w_t,\\
\hat x^u_0=&\ \mu_0, \quad \hat p_T=H\hat x^u_T+\bar H\IE x^u_T+L+\bar
L
\end{aligned}
\right.
\end{equation}
with
$$
\theta_{1,t}=\left(\bar a_t-b_tB^{-1}_t\bar D_t\right)\IE x^u_t-b_tB^{-1}_tG_t+\bar b_t,$$
$$
\begin{aligned}
\theta_{2,t}=&\ \left(\bar A_t-D^\top_tB^{-1}_t\bar D_t-\bar D^\top_tB^{-1}_tD_t\right.\\
             &\ \left.-\bar D^\top_tB^{-1}_t\bar D_t\right)\IE x^u_t+\left(\bar a^\top_t-\bar D^\top_tB^{-1}_tb^\top_t\right)\IE p_t\\
             &\ -\left(D_t+\bar D_t\right)^\top B^{-1}_tG_t+F_t+\bar F_t,
\end{aligned}
$$
where $\Sigma$ and $\bar{w}$ satisfy (\ref{sigma}) and (\ref{ip}),
and $\IE x^u$ and $\IE p$ solve (\ref{lhxm}) and (\ref{pxmr}),
respectively. We assume that the following condition holds.

\noindent \textbf{(A2).} There is a constant $C\geq0$ such that
$$D_t^\top B^{-1}_tD_t-A_t-CI_{n\times n}\leq0.$$

\noindent Then (\ref{lhsf}) has a unique solution $(\hat x^u, \hat p, \widehat{ \tilde{Q}})$ in the space $\mathscr L^2_{\mathscr F^{Y^u}}(0, T; \IR^{n+n+n\times\tilde r})$ by using Theorem 2.6 in Peng and Wu \cite{pw} again. Similarly,
let
\begin{equation}
\label{pxfr}\hat p_t=\Gamma_t\hat x^u_t+\Lambda_t
\end{equation}
for two deterministic and differentiable functions $\Gamma$ and
$\Lambda$ such that $\Gamma_T=H$ and $\Lambda_T=\bar H\IE
x^u_T+L+\bar L$. It follows from It\^o's formula that
$$
\begin{aligned}
d\hat p_t=&\ \dot\Gamma_t\hat x^u_tdt+\Gamma_td\hat x^u_t+\dot\Lambda_tdt\\
=&\ \left\{\left[\dot\Gamma_t+\Gamma_t\left(a_t-b_tB^{-1}_tD_t\right)-\Gamma_tb_tB^{-1}_tb^\top_t\Gamma_t\right]\hat x^u_t\right.\\
&\ \left.+\Gamma_t\left(\theta_{1,t}-b_tB^{-1}_tb^\top_t\Lambda_t\right)+\dot\Lambda_t\right\}dt\\
&\ +\Gamma_t\Sigma_tf_t^{\top}(h^{-1}_t)^{\top}d\bar w_t.
\end{aligned}
$$
Comparing it with the BSDE in (\ref{lhsf}), we derive
\begin{equation}
\label{gamma}\left\{\begin{aligned}
&\ \dot\Gamma_t+\Gamma_t\left(a_t-b_tB^{-1}_tD_t\right)+\left(a^\top_t-D^\top_tB^{-1}_tb^\top_t\right)\Gamma_t\\
&\ -\Gamma_tb_tB^{-1}_tb^\top_t\Gamma_t+A_t-D^\top_tB_tD_t=0,\\
&\ \Gamma_T=H
\end{aligned}
\right.
\end{equation}
and
\begin{equation}
\label{lambda}\left\{\begin{aligned}
&\ \dot\Lambda_t+\left(a^\top_t-D^\top_tB^{-1}_tb^\top_t-\Gamma_tb_tB^{-1}_tb^\top_t\right)\Lambda_t\\
&\ +\Gamma_t\theta_{1,t}+\theta_{2,t}=0,\\
&\ \Lambda_T=\bar H\IE x^u_T+L+\bar L,
\end{aligned}
\right.
\end{equation}
which have a unique solution, respectively. Substituting
(\ref{pxfr}) into (\ref{loc}), we get
\begin{equation}
\label{lofc}
u_t=-B^{-1}_t\left[\left(b^\top_t\Gamma_t+D_t\right)\hat x^u_t+\bar
D_t\IE x^u_t+b^\top_t\Lambda_t+G_t\right],
\end{equation}
where $\IE x^u$, $\Gamma$, $\Lambda$ and $\hat x^u$ solve
(\ref{lhxm}), (\ref{gamma}), (\ref{lambda}) and the closed-loop
system
\begin{equation}
\label{lcls}\left\{\begin{aligned}
d\hat x^u_t=&\ \left\{\left[a_t-b_tB^{-1}_t\left(D_t+b^\top_t\Gamma_t\right)\right]\hat x^u_t\right.\\
&\ \left. -b_tB^{-1}_tb^\top_t\Lambda_t+\theta_{1,t}\right\}dt+\Sigma_tf_t^{\top}(h^{-1}_t)^{\top}d\bar w_t,\\
\hat x^u_0=&\ \mu_0,
\end{aligned}
\right.
\end{equation}
respectively.

\emph{Step 4:} (\ref{lofc}) is the optimal control.

According to (\ref{ip}), $\bar w_t$ is an $\mathscr
F^{Y^0}_t$-adapted standard Brownian motion. Then it is easy to see
from (\ref{lcls}) that $\hat x^u_t$ is $\mathscr F^{Y^0}_t$-adapted,
and hence, $u_t$ given by (\ref{lofc}) is $\mathscr
F^{Y^0}_t$-adapted. On the other hand, applying It\^o's formula to
$\langle\hat x^u, \hat x^u\rangle$ with Burkholder-Davis-Gundy
inequality, we deduce $$\IE\sup_{0\leq t\leq T}|\hat
x^u_t|^2<+\infty.$$ Then $u\in\mathscr U^0_{ad}$. Next, we will
prove that $u_t$ is also $\mathscr F^{Y^u}_t$-adapted. If so, then
$u\in\mathscr U_{ad}$, and consequently, $u$ is optimal via Theorem
3.2. In fact, using (\ref{ip}) again, (\ref{lcls}) can be rewritten
as
$$
\left\{\begin{aligned}
d\hat x^u_t=&\ \left\{\left[a_t-b_tB^{-1}_t\left(D_t+b^\top_t\Gamma_t\right)\right]\hat x^u_t\right.\\
            &\ \left.-b_tB^{-1}_tb^\top_t\Lambda_t+\theta_{1,t}\right\}dt\\
            &\ +\Sigma_tf_t^{\top}(h^{-1}_t)^{\top}h^{-1}_t\left[dY^u_t\right.\\
            &\ \left.-\left(f_t\hat x^u_t+\bar f_t\IE x^u_t+g_t\right)dt\right],\\
\hat x^u_0=&\ \mu_0.
\end{aligned}
\right.
$$
From the optimal filtering equation, it is easy to check that $\hat
x^u_t$ is $\mathscr F^{Y^u}_t$-adapted, so is $u_t$. Then
$u\in\mathscr U_{ad}$. Therefore, the claim holds.

\emph{Step 5:} Optimal cost functional.

Since the solution $\Sigma$ of (\ref{sigma}) is independent of $v$, the optimal cost functional is rewritten as
$$
\begin{aligned}
J[u]=&\ J_1[u]+\int_0^TN^{\top}\chi_{0}^{t}\bar{\psi}_{t}dt+\frac{1}{2}\int_0^Ttr(A_t\Sigma_t)dt\\
     &\ +tr(H\Sigma_T)
\end{aligned}
$$
with
\begin{equation}
\label{lqcf2}\begin{aligned}
J_1[u]=&\ \frac{1}{2}\IE\left\{\int_0^T[\langle A_t\hat x^u_t, \hat x^u_t\rangle+\langle \bar A_t\IE x^u_t, \IE x^u_t\rangle\right.\\
&\ +\langle B_tu_t, u_t\rangle+2\langle D_t\hat x^u_t,u_t\rangle+2\langle\bar D_t\IE x^u_t, u_t\rangle\\
&\ +2\langle F_t,x^u_t\rangle+2\langle\bar F_t, \IE x^u_t\rangle+2\langle G_t, u_t\rangle]dt\\
&\ +\langle H\hat x^u_T, \hat x^u_T\rangle+\langle\bar H\IE x^u_T, \IE x^u_T\rangle+2\langle L, x^u_T\rangle\\
&\ \left.+2\langle\bar L, \IE x^u_T\rangle\right\}.
\end{aligned}
\end{equation}
Here $\hat x^u$ solves (\ref{lcls}), and $tr(A)$ denotes the
trace of the matrix $A$. Using It\^o's formula, we deduce
$$
\begin{aligned}
d\langle\Gamma_t\hat x^u_t, \hat x^u_t\rangle=&\ \left\langle\left(D^\top_tB_tD_t-A_t-\Gamma_tb_tB^{-1}_tb^\top_t\Gamma_t\right)\hat x^u_t\right.\\
&\ \left.-2\Gamma_tb_tB^{-1}_tb^\top_t\Lambda_t+2\Gamma_t\theta_{1,t}, \hat x^u_t\right\rangle dt\\
&\ +\left\langle\Sigma_tf_t^{\top}(h^{-1}_t)^{\top}, \Gamma_t\Sigma_tf_t^{\top}(h^{-1}_t)^{\top}\right\rangle dt\\
&\ +\left\langle\Gamma_t\Sigma_tf_t^{\top}(h^{-1}_t)^{\top}d\bar w_t, \hat x^u_t\right\rangle\\
&\ +\left\langle\Gamma_t\hat x^u_t, \Sigma_tf_t^{\top}(h^{-1}_t)^{\top}d\bar w_t\right\rangle.
\end{aligned}
$$
Then
\begin{equation}
\label{gammax}\begin{aligned}
&\ \IE\langle H\hat x^u_T, \hat x^u_T\rangle\\
=&\ \IE\int_0^T\left\langle\left(D^\top_tB_tD_t-A_t-\Gamma_tb_tB^{-1}_tb^\top_t\Gamma_t\right)\hat x^u_t\right.\\
&\ \left.-2\Gamma_tb_tB^{-1}_tb^\top_t\Lambda_t+2\Gamma_t\theta_{1,t}, \hat x^u_t\right\rangle dt\\
&\ +\int_0^T\left\langle\Sigma_tf_t^{\top}(h^{-1}_t)^{\top}, \Gamma_t\Sigma_tf_t^{\top}(h^{-1}_t)^{\top}\right\rangle dt\\
&\ +\langle\Gamma_0\mu_0, \mu_0\rangle.
\end{aligned}
\end{equation}
Similarly, applying It\^o's formula to $\langle\Lambda, \hat
x^u\rangle$, we get
\begin{equation}
\label{lambdax}\begin{aligned}
&\ \IE\langle\bar H\IE x^u_T+L+\bar L, \hat x^u_T\rangle\\
=&\ \int_0^T\left\langle \Lambda_t, \theta_{1,t}-b_tB^{-1}_tb^\top_t\Lambda_t\right\rangle dt+\langle\Lambda_0, e_0\rangle\\
&\ -\IE\int_0^T\langle\Gamma_t\theta_{1,t}+\theta_{2,t}, \hat
x^u_t\rangle dt.
\end{aligned}
\end{equation}
Plugging (\ref{lofc}), (\ref{gammax}) and (\ref{lambdax}) into
(\ref{lqcf2}), we derive
\begin{equation}
\label{lqcf3}\begin{aligned}
J[u]=&\ \frac{1}{2}\left\{\int_0^T\left\langle\left[\left(2D_t+\bar D_t\right)^\top B^{-1}_t\bar D_t-\bar A_t\right.\right.\right.\\
&\ \left.\left.\left.-2\bar a^\top_t\Gamma_t\right]\IE x^u_t, \IE x^u_t\right\rangle dt-\left\langle \bar H\IE x^u_T, \IE x^u_T\right\rangle\right\}\\
&\ +\int_0^T\left\langle\bar D^\top_tB^{-1}_tb^\top_t\Gamma_t-\Gamma_tb^\top_tB^{-1}_tG_t\right.\\
&\ \left.-\Gamma_tb_tB^{-1}_tb^\top_t\Lambda_t, \IE x^u_t\right\rangle dt\\
&\ +\int_0^T\left\langle b^\top_t\Gamma_t, B^{-1}_t\left(b^\top_t\Lambda_t+G_t\right)\right\rangle dt\\
&\ +\frac{1}{2}\int_0^T\left\langle \Lambda_t, 2\bar b_t-2b_tB^{-1}_tG_t-b_tB^{-1}_tb^\top_t\Lambda_t\right\rangle dt\\
&\ -\frac{1}{2}\int_0^T\langle G_t, B^{-1}_tG_t\rangle dt\\
&\ +\frac{1}{2}\langle\Gamma_0\mu_0, \mu_0\rangle+\langle\Lambda_0, \mu_0\rangle\\
&\ +\frac{1}{2}\int_0^T\left\langle\Sigma_tf_t^{\top}(h^{-1}_t)^{\top}, \Gamma_t\Sigma_tf_t^{\top}(h^{-1}_t)^{\top}\right\rangle dt\\
&\ +\int_0^TN^{\top}\chi_{0}^{t}\bar{\psi}_{t}dt+\frac{1}{2}\int_0^Ttr(A_t\Sigma_t)dt\\
&\ +tr(H\Sigma_T).
\end{aligned}
\end{equation}

We summarize the above deduction as follows.

\textbf{Proposition 4.2.} \emph{Under the assumptions (A1), (A2), $M=0$, and $\bar\beta_t=\gamma_t=\bar\gamma_t=\tilde\gamma_t=\bar{\tilde\gamma}_t=0$, the optimal feedback control and the optimal cost
functional of Problem (MFLQ) are explicitly given by (\ref{lofc}) and (\ref{lqcf3}), respectively.}

\section{Solution to Problem (AL)}

In this section, we are interested in explicitly computing Problem (AL) introduced in Section 1.2 with the assumption that the generator $$G(t, y, z)=\beta_ty+\psi_tv,$$ where $\beta$ and $\psi$ are deterministic and uniformly bounded. With the generator, Problem (AL) is a special case of Problem (MFLQ). (\ref{lhsm}) is reduced to a decoupled forward-backward ODE
$$
\left\{\begin{aligned}
\dot{\IE x^u_t}=&\ (a_t+\bar a_t)\IE x^u_t-B^{-1}_tb^2_t\IE p_t\\
&\ +NB^{-1}_tb_t\psi_te^{\int_0^t\beta_sds}+\bar b_t,\\
\dot{\IE p_t}=&\ -(a_t+\bar a_t)\IE p_t,\\
\IE x^u_0=&\ \mu_0, \quad \IE p_T=-Ne^{\int_0^T\beta_tdt}.
\end{aligned}
\right.
$$
Solving it, we get
\begin{equation}
\label{al:p}
\IE p_t=-Ne^{\int_0^T\beta_tdt+\int_t^T(a_s+\bar a_s)ds}
\end{equation}
and
\begin{equation}
\label{al:x}
\begin{aligned}
\IE x^u_t=&\ \mu_0e^{\int_0^t(a_s+\bar a_s)ds}-\int_0^tB^{-1}_sb^2_se^{\int_s^t(a_r+\bar a_r)dr}\IE p_sds\\
&\ +\int_0^t\left(\bar b_s+NB^{-1}_sb_s\psi_se^{\int_0^s\beta_rdr}\right)e^{\int_s^t(a_r+\bar a_r)dr}ds,
\end{aligned}
\end{equation}
where $\IE p$ is determined by (\ref{al:p}). (\ref{gamma}) and (\ref{lambda}) are rewritten as
\begin{equation}
\label{al:gamma}\left\{\begin{aligned}
&\ \dot\Gamma_t+2a_t\Gamma_t-B^{-1}_tb^2_t\Gamma^2_t=0,\\
&\ \Gamma_T=H
\end{aligned}
\right.
\end{equation}
and
\begin{equation}
\label{al:lambda}\left\{\begin{aligned}
&\ \dot\Lambda_t+\left(a_t-B^{-1}_tb^2_t\Gamma_t\right)\Lambda_t+\Gamma_t\theta_{1,t}+\theta_{2,t}=0,\\
&\ \Lambda_T=-H\IE x^u_T-Ne^{\int_0^T\beta_tdt},
\end{aligned}
\right.
\end{equation}
where
$$
\theta_{1,t}=\bar a_t\IE x^u_t+B^{-1}_tb_t\psi_tNe^{\int_0^t\beta_sds}+\bar b_t, \quad\theta_{2,t}=\bar a_t\IE p_t,
$$
with $\IE p$ and $\IE x$ be determined by (\ref{al:p}) and (\ref{al:x}). This gives
$$
\begin{aligned}
\Lambda_t=&\ -\left(H\IE x^u_T+Ne^{\int_0^T\beta_sds}\right)e^{\int_t^T(a_s-B^{-1}_sb^2_s\Gamma_s)ds}\\
&\ +\int_t^T(\Gamma_s\theta_{1,s}+\theta_{2,s})e^{\int_t^s(a_r-B^{-1}_rb^2_r\Gamma_r)dr}ds.
\end{aligned}
$$
Then the optimal control strategy in a feedback form is
$$
u_t=-B^{-1}_t\left[b_t\left(\Gamma_t\hat x^u_t+\Lambda_t\right)-N\psi_te^{\int_0^t\beta_sds}\right].
$$
Here $\Gamma$ and $\Lambda$ solve (\ref{al:gamma}) and (\ref{al:lambda}). The filtered cash-balance $\hat x^u$ satisfies a closed-loop system
$$
\left\{\begin{aligned}
d\hat x^u_t=&\ \left[\left(a_t-B^{-1}_tb^2_t\Gamma_t\right)\hat x^u_t-B^{-1}_tb^2_t\Lambda_t+\theta_{1,t}\right]dt\\
            &\ +\Sigma_tf_th^{-1}_td\bar w_t,\\
\hat x^u_0=&\ \mu_0
\end{aligned}
\right.
$$
with $\Sigma$ being governed by
$$
\left\{\begin{aligned}
&\ \dot\Sigma_t-2a_t\Sigma_t+h^{-2}_tf^2_t\Sigma^2_t-c^2_t=0,\\
&\ \Sigma_0=\sigma_0.
\end{aligned}
\right.
$$

The remaining part of this section is devoted to illustrating the above results via numerical computations. For the sake of illustrations, let the time granularity for all the market parameters be yearly. Suppose that the initial investment $\xi=1$, the discount rate $a_t=0.03$, the liability rate $\bar b_t=0.01$, and the volatility rates $c_t=0.04$ and $h_t=0.1$. We further suppose that the market parameters $\bar a_t=0.03$, $\beta_t=0.06$, $f_t=0.1$, $b_t=\psi_t=B_t=N=1$, and $H=0.01$. With the data, it is easy to see that
$$
\begin{aligned}
\IE p_t=&\ -e^{0.06(2-t)}, \\
\IE x^u_t=&\ e^{0.06t}\left[1+t+\frac{25}{3}e^{0.12}(1-e^{-0.12t})\right.\\
&\ \left.+\frac{1}{6}(1-e^{-0.06t})\right],\\
\theta_{1,t}=&\ 0.03\IE x^u_t+e^{0.06t}+0.01, \quad \theta_{2,t}=-0.03e^{0.06(2-t)},\\
\Gamma_t=&\ \frac{0.06e^{0.06(1-t)}}{5+e^{0.06(1-t)}},
\end{aligned}
$$
and
$$
\begin{aligned}
\Lambda_t=&\ -(0.01\IE x^u_1+e^{0.06})e^{\int_t^1(0.03-\Gamma_s)ds}\\
          &\ +\int_t^1(\Gamma_s\theta_{1,s}+\theta_{2,s})e^{\int_t^s(0.03-\Gamma_r)dr}ds.
\end{aligned}
$$
Then the optimal control strategy of the firm is
$$
u_t=-(\Gamma_t\hat x^u_t+\Lambda_t-e^{0.03t}),
$$
where the optimal filtering of the cash-balance $x^u$ satisfies
$$\left\{\begin{aligned}
d\hat x^u_t=&\ \left[\left(0.03-\Gamma_t\right)\hat x^u_t-\Lambda_t+\theta_{1,t}\right]dt+\Sigma_td\bar w_t,\\
\hat x^u_0=&\ 1
\end{aligned}
\right.
$$
with $$\Sigma_t=\frac{0.08(e^{0.1t}-1)}{e^{0.1t}-4},$$ $\Gamma$, $\Lambda$ and $\theta_1$ being determined above.

\section{Concluding remarks}

This article studies an optimal control problem for mean-field FBSDE with noisy observation. Since mean-field FBSDE and optimal filtering are considered, the control problem has been basically unexplored so far. The control problem covers more models in reality, but causes a trouble in solving the problem. The backward separation method with the decomposition of the state and the observation is further developed, and is introduced to overcome the resulting difficulty. These results obtained in this article improve the first author's previous works \cite{wwx1,wwx2,wwz1,wzz}, and are helpful for studying mean-field game for FBSDE and systematic risk model with noisy observation. The details of how we study these problems will be presented elsewhere.

Let us now make several remarks in order to close this section. (1) In most optimal control problems for mean-field stochastic systems, we assume that all coefficients of the optimal control problems are deterministic. Otherwise, there is an immediate difficulty to study the problems. One reason is that the key equality $\IE(a_tx_t)=a_t\IE x_t$ is no longer true if $a_t$ is also a stochastic process. But some special cases with stochastic coefficients, say, Example 2.1, can be solved by a simple reduction method. Then it is natural to ask if the method is applicable for slightly more complicated cases. We hope to answer it in the near future. (2) Similar to Example 2.2, Problem (MFLQ) can also be reduced to an LQ problem for non-standard FBSDE with control set constraint. This motivates us to investigate such a class of LQ problems for non-standard FBSDEs in the future. In return, it will be helpful to study LQ problems for mean-field FBSDEs. (3) The solution of the BSDE in (\ref{se}) is a non-Gaussian process in general, and thus, the optimal filter of the BSDE is infinite. Then it is highly desirable to study the numerical approximation of the optimal filter and the optimal control in our future publications.

\begin{ack}                               
The authors would like to thank the editors and three anonymous referees for their constructive and insightful comments for improving the quality of
this work.
\end{ack}

\section*{Appendix}

We present three lemmas first, and then give a proof of Theorem 2.1.

\textbf{Lemma A.1.}  \emph{For any $v_j\in\mathscr L^2_{\mathscr F}(0, T; \IR^k)$, $j=1, 2$, there is a constant
$C>0$ such that}
$$\IE\sup_{0\leq t\leq T}|x^{v_1}_t-x^{v_2}_t|^2\leq C\IE\int_0^T|v_{1,t}-v_{2,t}|^2dt,$$
\begin{equation*}
\begin{aligned}
       & \IE\sup_{0\leq t\leq T}|y^{v_1}_t-y^{v_2}_t|^2\\
       & +\IE\int_{0}^{T}\left(|z^{v_1}_t-z^{v_2}_t|^2+|\tilde{z}^{v_1}_t-\tilde{z}^{v_2}_t|^2\right)dt\\
  \leq & C\left( \IE|x^{v_1}_T-x^{v_2}_T|^2+\int_{0}^{T}\sup_{0\leq s\leq t}\IE|x^{v_1}_s-x^{v_2}_s|^2dt\right.\\
  & \left.+\IE\int_0^T|v_{1,t}-v_{2,t}|^2dt\right).
\end{aligned}
\end{equation*}

\textbf{Proof.} These two estimates can be derived by It\^o's formula, Gronwall's inequality and
Burkholder-Davis-Gundy inequality. We omit the proof for simplicity. \hfill$\Box$

\textbf{Lemma A.2.} \emph{For any $v$, $v_j\in\mathscr U_{ad}$ $(j=1, 2, \cdots)$ satisfying
$v_j\rightarrow v$ in $\mathscr L^2_{\mathscr F}(0, T; \IR^k)$, it holds}
$$\lim_{j\rightarrow+\infty}J[v_j]=J[v].$$

\textbf{Proof.} Using Taylor's expansion, H\"{o}lder's
inequality and Lemma A.1, we deduce
$$
\begin{aligned}
&\left|\IE\int_0^Tl(\Theta^{v_j}_t)dt-\IE\int_0^Tl(\Theta^v_t)dt\right|\\
\leq& C\IE\int_0^T\left(1+|x^{v_j}_t|+|x^v_t|+\IE|x^{v_j}_t|+\IE|x^v_t|\right.\\
&\ \left. +|v_{j,t}|+|v_t|\right)\left(|x^{v_j}_t-x^v_t|+|\IE x^{v_j}_t-\IE x^v_t|\right.\\
&\ \left. +|v_{j,t}-v_t|\right)dt\\
\leq& C\sqrt{\IE\int_0^T\aleph^{v, v_j}_tdt}\left(\sqrt{\IE\sup_{0\leq t\leq T}|x^{v_j}_t-x^v_t|^2}\right.\\
&\ \left.+\sqrt{\IE\int_0^T|v_{j,t}-v_t|^2dt}\right)\rightarrow0
\end{aligned}
$$
as $j\rightarrow+\infty$, where $C>0$ is a constant, and
$$
\begin{aligned}
\aleph^{v, v_j}_t=&\ 1+|x^{v_j}_t|^2+|x^v_t|^2+\IE|x^{v_j}_t|^2\\
&\ +\IE|x^v_t|^2+|v_{j,t}|^2+|v_t|^2.
\end{aligned}
$$

In a same way, we have $$\IE\phi(\Xi^{v_j}_T)\rightarrow\IE\phi(\Xi^v_T),\qquad \IE\varphi(y^{v_j}_0)\rightarrow\IE\varphi(y^v_0)$$ with $j\rightarrow+\infty$. Then the proof is complete.       \hfill$\Box$

\textbf{Lemma A.3.} \emph{$\mathscr U_{ad}$ is dense in $\mathscr
U^0_{ad}$.}

\textbf{Proof.} For any $v\in\mathscr U^0_{ad}$, define a family of
controls by
\[v_{j\,t}=\left\{
\begin{array}{l l}
\nu, & \hbox{for}\ 0\leq t\leq\delta_j,\\
\frac{1}{\delta_j}\int_{(i-1)\delta_j}^{i\delta_j}v_sds,           & \hbox{for}\ i\delta_j<t\leq(i+1)\delta_j,
\end{array}
\right.
\]
where $\nu\in U$, $i$, $j$ are natural numbers, $1\leq i\leq j-1$,
and $\delta_j=T/j$. Similar to Bensoussan \cite{ba1}, we can prove that (i) $v_j\in\mathscr U_{ad}$ for any $j$, and (ii)
$v_j\rightarrow v$ as $j\rightarrow+\infty$ in
$\mathscr L^2_{\mathscr F^{Y^0}}(0, T; U)$. Then it implies the desired result.
\hfill$\Box$

\textbf{Proof of Theorem 2.1.} From Definition 2.1, we have $\mathscr U_{ad}\subseteq\mathscr U^0_{ad}$, and thus, $\inf_{v'\in\mathscr U_{ad}}J[v']\geq\inf_{v\in\mathscr U^0_{ad}}J[v].$
On the other hand, since $v_j$ defined in the proof of Lemma A.3 is an element of $\mathscr U_{ad}$, then $\inf_{v'\in\mathscr U_{ad}}J[v']\leq J[v_j]$, and consequently, it follows from Lemma A.2 that $\inf_{v'\in\mathscr U_{ad}}J[v']\leq\lim_{j\rightarrow+\infty}J[v_j]=J[v]$. Due to the arbitrariness of $v$, then $\inf_{v'\in\mathscr U_{ad}}J[v']\leq\inf_{v\in\mathscr U^0_{ad}}J[v].$ Thus, the proof is complete. \hfill$\Box$


\end{document}